%% file: main.tex
\documentclass[pdflatex,sn-basic]{sn-jnl}
\setcitestyle{aysep={,}}

\usepackage{
    amssymb,
    amsmath,
    commath,
    cleveref,
    dsfont,
    subcaption,
    booktabs,
    multirow,
    graphicx,
    algorithm,
    algpseudocodex,
}

\usepackage[
    group-separator={,}
]{siunitx}

\input{tikz_preamble}

\title{Integrated Airline Fleet and Crew Recovery through Local Search}

\author*[1,2]{\fnm{Philip} \sur{de Bruin}}\email{p.debruin@uu.nl}% orcid: 0000-0002-1981-0527
\author[1,2]{\fnm{Marjan} \sur{van den Akker}}\email{j.m.vandenakker@uu.nl}% orcid: 0000-0002-7114-0655
\author[2]{\fnm{Kunal} \sur{Kumar}}\nomail
\author[1,2]{\fnm{Lisanne} \sur{Heuseveldt}}\nomail% orcid: 0009-0002-3879-6166

\affil[1]{\orgdiv{Department of Information and Computing Sciences}, \orgname{Utrecht University}, \city{Utrecht}, \country{The Netherlands}}
\affil[2]{\orgname{KLM Royal Dutch Airlines}, \city{Amstelveen}, \country{The Netherlands}}

\abstract{
\input{0_abstract.tex}
}


\begin{document}

\maketitle

\section{Introduction}\label{sec:introduction}
\input{1_introduction.tex}

\section{Literature Overview}\label{sec:literature}
\input{2_literature.tex}

\section{Problem Description}\label{sec:problem_description}
\input{3_problem_description.tex}

\section{Local Search Algorithm}\label{sec:local_search}
\input{4_local_search.tex}

\section{Experimental Results}\label{sec:results}
\input{5_results.tex}

\section{Conclusion}\label{sec:conclusion}
\input{6_conclusion.tex}

\section*{Acknowledgements}
We thank Marc Paelinck for sharing his knowledge and insights on crew scheduling.

This research was conducted as part of the ORDERbyCHAOS project (Optimizing Resilience, Mitigating Disruptions, and Enhancing Robustness by using Combinatorial Heuristics in Airline Operational Scheduling), which is a collaboration of KLM Royal Dutch Airlines and Utrecht University, and carried out within the KickStartAI program.

\appendix
\section{MIP Formulation for Tail Recovery}\label{sec:tail_mip}
\input{7_tail_mip.tex}

\bibliography{papers}

\end{document}

%% file: tikz_preamble.tex
\usepackage{xcolor}
\definecolor{klmlightblue}{HTML}{00A1DE}
\definecolor{klmdarkblue}{HTML}{003145}

\definecolor{uuYellow}{HTML}{FFCD00}
\definecolor{uuRed}{HTML}{C00A35}
\definecolor{uuCreme}{HTML}{FFE6AB}
\definecolor{uuOrange}{HTML}{F3965E}
\definecolor{uuBordeauxRed}{HTML}{AA1555}
\definecolor{uuBrown}{HTML}{6E3B23}
\definecolor{uuGreen}{HTML}{24A793}
\definecolor{uuBlue}{HTML}{5287C6}
\definecolor{uuDarkBlue}{HTML}{001240}
\definecolor{uuPurple}{HTML}{5B2182}

\usepackage{tikz, fontawesome5}
\usetikzlibrary{
    backgrounds,
    calc,
    decorations.pathreplacing,
    decorations.pathmorphing,
    patterns,
}
\tikzset{
    flight/.style = {
        rectangle, very thin, draw = black, fill = white,
    },
    outbound/.style = {
        outboundflight
    },
    inbound/.style = {
        inboundflight
    },
    selected/.style = {
        draw=uuOrange, very thick
    },
}
\tikzset{
    invisible/.style={opacity=0},
    visible on/.style={alt=#1{}{invisible}},
    alt/.code args={<#1>#2#3}{%
        \alt<#1>{\pgfkeysalso{#2}}{\pgfkeysalso{#3}} %
    },
}

\makeatletter

\pgfdeclareshape{inboundflight}{%
    \inheritsavedanchors[from=rectangle] %
    \inheritanchorborder[from=rectangle]
    \inheritanchor[from=rectangle]{center}
    \foreach \anchor in {north,north west,north east,center,west,east,mid, mid west,mid east,base,base west,base east,south,south west,south east,text}{%
        \inheritanchor[from=rectangle]{\anchor}}%
    \backgroundpath{%
        \northeast \pgf@xb=\pgf@x \pgf@ya=\pgf@y
        \southwest \pgf@xa=\pgf@x \pgf@yb=\pgf@y
        \pgf@xc=\pgf@xa \advance\pgf@xc by 3pt
        \pgfpathmoveto{\pgfpoint{\pgf@xa}{\pgf@ya}}
        \pgfpathlineto{\pgfpoint{\pgf@xb}{\pgf@ya}}
        \pgfpathlineto{\pgfpoint{\pgf@xb}{\pgf@yb}}
        \pgfpathlineto{\pgfpoint{\pgf@xa}{\pgf@yb}}
        \pgfpathlineto{\pgfpoint{\pgf@xa}{\pgf@yb + 0.40\ht\pgfnodeparttextbox}}
        \pgfpathlineto{\pgfpoint{\pgf@xc}{\pgf@yb + 0.40\ht\pgfnodeparttextbox}}
        \pgfpathlineto{\pgfpoint{\pgf@xc}{\pgf@yb + 0.80\ht\pgfnodeparttextbox}}
        \pgfpathlineto{\pgfpoint{\pgf@xa}{\pgf@yb + 0.80\ht\pgfnodeparttextbox}}
        \pgfpathlineto{\pgfpoint{\pgf@xa}{\pgf@yb + 1.20\ht\pgfnodeparttextbox}}
        \pgfpathlineto{\pgfpoint{\pgf@xc}{\pgf@yb + 1.20\ht\pgfnodeparttextbox}}
        \pgfpathlineto{\pgfpoint{\pgf@xc}{\pgf@yb + 1.60\ht\pgfnodeparttextbox}}
        \pgfpathlineto{\pgfpoint{\pgf@xa}{\pgf@yb + 1.60\ht\pgfnodeparttextbox}}
        \pgfpathclose
    }
}

\pgfdeclareshape{outboundflight}{%
    \inheritsavedanchors[from=rectangle] %
    \inheritanchorborder[from=rectangle]
    \inheritanchor[from=rectangle]{center}
    \foreach \anchor in {north,north west,north east,center,west,east,mid, mid west,mid east,base,base west,base east,south,south west,south east,text}{%
        \inheritanchor[from=rectangle]{\anchor}}%
    \backgroundpath{%
        \northeast \pgf@xa=\pgf@x \pgf@ya=\pgf@y
        \southwest \pgf@xb=\pgf@x \pgf@yb=\pgf@y
        \pgf@xc=\pgf@xa \advance\pgf@xc by 3pt
        \pgfpathmoveto{\pgfpoint{\pgf@xa}{\pgf@ya}}
        \pgfpathlineto{\pgfpoint{\pgf@xb}{\pgf@ya}}
        \pgfpathlineto{\pgfpoint{\pgf@xb}{\pgf@yb}}
        \pgfpathlineto{\pgfpoint{\pgf@xa}{\pgf@yb}}
        \pgfpathlineto{\pgfpoint{\pgf@xa}{\pgf@yb + 0.40\ht\pgfnodeparttextbox}}
        \pgfpathlineto{\pgfpoint{\pgf@xc}{\pgf@yb + 0.40\ht\pgfnodeparttextbox}}
        \pgfpathlineto{\pgfpoint{\pgf@xc}{\pgf@yb + 0.80\ht\pgfnodeparttextbox}}
        \pgfpathlineto{\pgfpoint{\pgf@xa}{\pgf@yb + 0.80\ht\pgfnodeparttextbox}}
        \pgfpathlineto{\pgfpoint{\pgf@xa}{\pgf@yb + 1.20\ht\pgfnodeparttextbox}}
        \pgfpathlineto{\pgfpoint{\pgf@xc}{\pgf@yb + 1.20\ht\pgfnodeparttextbox}}
        \pgfpathlineto{\pgfpoint{\pgf@xc}{\pgf@yb + 1.60\ht\pgfnodeparttextbox}}
        \pgfpathlineto{\pgfpoint{\pgf@xa}{\pgf@yb + 1.60\ht\pgfnodeparttextbox}}
        \pgfpathclose
    }
}

\makeatother

%% file: 0_abstract.tex
Airline operations are prone to delays and disruptions, since the schedules are generally tight and depend on many resources.
When disruptions occur, the flight schedule needs to be adjusted such that the operation can continue.
Since this happens during the day of operations, this needs to be done as close to real time as possible, posing a challenge with respect to computation time.
Moreover, to limit the impact of disruptions, we want a solution with minimal cost and passenger impact.
Since airline operations include many interlinked decisions, an integrated approach leads to better overall solutions.
We specifically look at resolving these disruptions in both the aircraft and crew schedules.
Resolving these disruptions is complex, especially when it is done in an integrated way, i.e. including multiple different resources.

To solve this problem in an integrated manner, we developed a fast simulated annealing approach.
This approach is compared with traditional approaches, and an experimental study is done to evaluate different neighbour generation methods, and to investigate different recovery scenarios and strategies.
The comparison is done using real-world data from KLM Royal Dutch Airlines.
Here, we show that our approach resolves disruptions quickly and in a cost-efficient manner, and that it outperforms traditional approaches.
Compared to naive delay propagation, our method saves 40\% in non-performance costs.
Moreover, while most airlines use tools that consider resources separately, our approach shows that integrated disruption management is possible within 30 seconds.

%% file: 1_introduction.tex
Airline operations are prone to disruptions.
These could be due to bad weather conditions, airport congestion, or other unforeseen circumstances.
Airline schedules are generally tight and depend on many resources. Thus disruptions can have substantial effects because they propagate through the schedule.
Improper handling of disruptions can lead to many cancellations, which can cause even more chaos and large recovery costs.
This creates additional operational challenges as it becomes harder to resume the original schedule.
Disruptions also have a big impact on passengers, who need to be reaccommodated.
In the EU, this translates to an estimated cost of 6.5 billion euros in 2024 \citep{Skycop2025}.
On top of that, the negative impact on passenger satisfaction creates a competitive disadvantage for the airline.
Overall, it is estimated that disruptions cost airlines around 8\% of their total revenue \citep{CMAC2024}.
Therefore, it is important to handle these disruptions correctly and efficiently, as they can become very costly.
Since the COVID pandemic, most airlines have been left with a lower availability of their resources, making airline operations even more vulnerable to disruptions.
This is clearly visible in a recent performance review by Eurocontrol \citep{PRREuroControl2024}.
They show that punctuality has been on a decline since the summer of 2021 and punctuality is also lower than 2019, while overall traffic levels are below 2019 levels.
Examples of resource issues are: shortages of flight crew and supply chain issues on aircraft parts, where the latter can cause aircraft to stay in maintenance longer and they are thus unavailable for flight.

When a disruption occurs, an airline typically has a range of options to handle the situation. 
They can make changes to the resource assignment, essentially moving or swapping flights between, for example, aircraft.
For example, when maintenance of an aircraft is delayed, it may not be able to perform the first flight planned after maintenance.
In this situation, the airline may choose to reassign the flight to another aircraft, or swap assignments between two aircraft.
Furthermore, they can make adjustments to the schedule, i.e. to delay or, in an extreme case, cancel the flight.
These options each come with their own impact, including on passengers, schedule stability, and other key considerations, all of which can be represented in the recovery objective through associated costs.
To avoid large losses, besides resolving a disruption, we also need to minimize these costs.
Next to that, disruptions happen live during the day of operation and may be revealed at short notice.
Consequently, they need to be resolved quickly to avoid grounding operations; and any changes made to the schedule need to be communicated to not only the passengers but also other operational stakeholders, such as ground operations, airport, air traffic control, and others, as they need to adjust their schedules to incorporate these changes.
The need to communicate changes to all these operational stakeholders adds a lot of time pressure for those involved in resolving these disruptions.
For this reason, a schedule in which all disruptions are resolved needs to be found within at most a few minutes.

This is known as the \textsc{Airline Recovery Problem} and is usually divided into aircraft, crew, and passenger recovery, each focussing on their respective part of the process.
Since flight delays and cancellations are changes to the flight schedule, the inclusion of those actions is often referred to as schedule recovery.

To solve the \textsc{Airline Recovery Problem}, one can consider separate recovery problems for each of the resources and solve them sequentially.
Here, each part boils down to a hard optimization problem with its own computational requirements.
However, these problems have a lot of interlinked decisions, meaning that decisions in one stage may lead to conflicts or infeasibilities in the other.
For example, the delay of a flight might be preferable for the aircraft assignment, but may mean that the assigned crew exceeds their maximum allowed duty hours.
Considering an integrated approach mitigates these issues, resulting in better overall solutions.
However, integrating decision-making with multiple resources into one model adds more complexity, increasing the model size considerably, which makes it more challenging to find a good solution in a short amount of time.
Although many papers focus on a single resource, the integrated approach has been studied before by, for example, \citet{Petersen2012}, and \citet{Arikan2017}.
Here, the disruptions under consideration are mainly related to the reduction of airport capacity.
There has hardly been any attention to the unexpected unavailability of aircraft or crew members.
These studies also show models with a typical runtime up to  15 to 20 minutes, meaning that application in practice is limited.
Recent machine learning approaches, such as \citet{Ding2023}, close this runtime performance gap, with a typical runtime of several seconds.

Earlier works show limitations with respect to crew duty regulations, mainly regarding the rest periods.
Either the rest period is seen as a fixed block, or the recovery algorithm does not consider a sufficiently long planning horizon to account for rest periods.
However, according to European regulation, the time required for the rest period has to depend on the duty length preceding this rest.
Lastly, in most of the earlier work, only a single aircraft type is considered, meaning that all crew can fly on all the aeroplanes available.
When swapping aircraft of different types, however, we have to deal with complications due to crew not being licensed to fly a certain aircraft type.
Thus, changing the aircraft type also means that the right crew swaps need to happen, increasing the complexity of the problem.

Our contributions are twofold.
First, we develop a simulated annealing algorithm for the {\sc Airline Recovery Problem} that jointly considers aircraft schedules, crew schedules, and passenger related costs, including missed connections.
The algorithm has an average runtime below 30 seconds, making it suitable for time-critical operational use.
Second, our algorithm handles the three operational aspects that received limited attention in earlier work.
It accounts for resource unavailability and incorporates European regulations for rest and duty times, whereas previous work mainly considers fixed blocks for duty and rest times.
Finally, our algorithm allows swaps of aircraft of different types.
This recovery option is difficult to use in a sequential approach because it affects both the crew and aircraft schedules.

The research is done in close collaboration with KLM Royal Dutch Airlines.
We perform an elaborate experimental study based on real-life data from this airline.
We investigate different disruptions and compare different neighbour generation methods.
Furthermore, we investigate the effect of penalizing changes in the schedule and the use of different delay cost models.
Lastly, we compare our algorithm with different sequential approaches.
This includes the use of a MIP to resolve disruptions in the aircraft schedule, an approach that is currently often done in practice.

The remainder of the paper is organized as follows.
We first give an overview of the current literature on this topic in \Cref{sec:literature}.
Then, we discuss the full problem in \Cref{sec:problem_description}.
In \Cref{sec:local_search}, we go into the details of the local search we developed.
Lastly, we provide a computational study in \Cref{sec:results}, and the conclusion in \Cref{sec:conclusion}.

%% file: 2_literature.tex
In this section, we give an overview of the literature related to \textsc{Airline Recovery Problem}.
Recently, \citet{Wu2025} and \citet{Hassan2021} presented surveys on this problem.
\Citet{Hassan2021} also discussed the practical challenges in this area.
They both show that in recent years there has been an increasing focus on integrated recovery.

Since aircraft were traditionally the most scarce resource of an airline, there has been a lot of research on aircraft and schedule recovery. Furthermore, the use of these aircraft is expensive, so it is important to use them as efficiently as possible.
\Citet{Aguiar2011} compared three different algorithms, namely a hill-climbing algorithm, a simulated annealing algorithm, and a genetic algorithm.
These algorithms are based on swapping the assignments of two aircraft.
They find that the genetic algorithm performs best;
compared to simulated annealing, it finds slightly better solutions quicker.
\Citet{Guimarans2015} use a large neighbourhood search.
In their search procedure, they remove the assignments of specific aircraft and then rebuild a solution.
Both these papers consider disruptions due to factors outside the airline's control.
Thus, delays caused by, for example, weather conditions or airport congestion.

\Citet{Zhu2015} create a stochastic recovery model, where they consider the aircraft maintenance as a potential disruption source.
They then solve the problem in a recourse model.
Here, the here-and-now decision is to assign aircraft to flights with the option to cancel flights, and the wait-and-see decision delays flights.
\Citet{Lee2020} combine robust planning with recovery from disruptions.
They use a queuing model in order to predict delays at hub airports in order to create a robust flight schedule.
This can be combined with an ILP in order to react to disruptions.

\Citet{Rashedi2025} solve the aircraft recovery problem by employing a machine learning model to fix parts of the flight network.
This prunes the solution space, after which the remaining network is solved using an MIP.
Using this approach, they obtain near-optimal solutions within approximately five minutes.

Integrating aircraft recovery with crew recovery is considerably more complex.
Not only do we need to make more decisions, but crew recovery also introduces additional complex constraints due to flying time regulations.
\Citet{Petersen2012} study the fully integrated recovery problem.
Thus, they consider fleet, crew, and passengers.
For this, they develop a complicated algorithm making use of both Bender's cuts and column generation.
Their model only optimizes passenger recovery, and does not solve for global optimality.
Furthermore, their recovery window ends at the current day limiting crew recovery options due to not being able to consider their rest period.
They test their model on a real-world network containing about \num{800} daily flights and achieve a runtime of about $30$ minutes.

\Citet{Brunner2014} studies the effect of the Ground Delay Program (GDP) in the US, which can be viewed as a specific type of disruption in the sense that the Ground Delay Program may change the available slots at an airport. This means that flights need to be rescheduled to the new slots.
For this, \citet{Brunner2014} built an MIP model considering passenger and crew connections.
In their experimental results, they study a single GDP containing 71 flights.
This GDP is then solved within a few seconds.

\Citet{Zhang2015} look at fleet and crew recovery but do not consider passengers.
They solve an MILP for aircraft recovery with crew considerations and an MILP for crew recovery with aircraft considerations.
They test solving these models both in a sequential and iterative manner.
Their method resolves disruptions on a real-world flight network of roughly \num{350} flights within $2$ minutes.
In later work, \citet{Maher2016} uses a column-and-row generation approach to solve the integrated crew and fleet recovery problem.
They test their approach on generated instances containing about \num{440} flights.
Their method solves most of these within the time limit of $20$ minutes.

\Citet{Arikan2017} are the most complete when it comes to integrated recovery.
They consider fleet, crew, and passengers, and they have more recovery actions compared to other papers.
However, their network is less flexible with respect to crew rest periods.
Since they allow for different cruising speeds, they need a conic quadratic mixed integer formulation.
In order to reduce the overall size of the model, they develop a reduction algorithm to eliminate recovery actions that are unlikely to be taken.
Without the reduction algorithm, they can solve disruptions in a flight network containing \num{473} flights in about $18$ minutes, where $93\%$ of these instances were solved to optimality.
However, when using a flight network containing \num{1254} flights, they had to use this reduction algorithm as the instance became too big.
In that case, they could solve disruptions in about $15$ minutes.
\Citet{Xu2023} adapt the MIP used by \citet{Arikan2017} in order to account for disruptions due to disease spreading.

Approaches taking full consideration of (European) crew rest and duty regulations are very new.
\Citet{Davies2026} show a model on integrated recovery used by the Lufthansa Group.
They adapt column generation to optimally combine solutions from the different subproblems (i.e. the tail, crew, and passenger recovery problems).
These subproblems are then solved using heuristic and exact methods.
Their models are used in practice, but they don't report how fast certain instance sizes can be solved.
\Citet{Puchkova2026} use logic-based row generation in order to solve fleet, crew, and passenger recovery.
With this approach, they solve disruptions with up to $467$ flights in $79$ minutes.
Both of these works closely follow the same crew rest and duty regulations as we do.

In recent years, approaches that include machine learning have been introduced.
\Citet{Eikelenboom2023} use machine learning to reduce the problem size. They use a ranking approach to learn which parts of the network to consider when resolving disruptions.
This vastly reduced runtime, resolving disruptions on flight networks consisting of about $327$ flights in $70$ seconds.
Their solutions seem to contain a relatively large number of cancellations.
\Citet{Ding2023} also use a machine learning approach to perform airline recovery.
Their approach is quite different; they use reinforcement learning with a variable neighbourhood search.
Then, using proximal policy optimization, the model learns which neighbourhood operations to select in a given state.
Their recovery model is similar to that of \citet{Arikan2017}, allowing the same recovery actions.
They test their method on synthetic data, where their biggest instance of \num{931} flights, is solved in $15$ seconds.
\Citet{Wang2025} also build a deep reinforcement learning approach based upon the flight network of \citet{Arikan2017}.
Their method then uses the same synthetic data generation as \citet{Ding2023}, and solves instances with up to \num{30} aircraft in around $10$ seconds.

As mentioned in \Cref{sec:introduction}, the disruptions covered in these works are focused around a reduction in airport capacity.
Thus, there is not much focus on unavailability of aircraft or crew.
Next to that the modelling of crew regulations remains limited, mainly focussing on maximum duty times.
Recent work has begun to address these limitations.
Furthermore, except for \citet{Petersen2012}, multiple aircraft types are not considered.

%% file: 3_problem_description.tex
In the \textsc{Airline Recovery Problem}, there is an initial assignment of the flights to the resources, i.e. aircraft and crew.
However, this assignment is currently infeasible due to disruptions.
We need to make it feasible again, which can be done by making certain changes.
These changes, so-called `levers', will be explained later.
Each of these changes influences the cost of the schedule, be it either operational costs or costs related to passengers.
The goal is to recover the given schedule, where we want to find a feasible schedule with minimal cost.

We consider a schedule for a set of flights $F$.
Each flight $f \in F$ has a departure time $f_d$ and an arrival time $f_a$, and is performed by a specific aircraft (tail) and a set of crew members.
Let $T$ and $C$ respectively be the sets of tails and crew members available.
The set $T$ contains aircraft of different types and with different passenger capacities.
A crew member in $C$ is only licensed to fly with aircraft of certain types and therefore with a subset of tails in $T$.
In our case, we only consider cockpit crew, thus every $s \in C$ represents one of the pilots.
Note that our approach is generalizable to cabin crew as well.
The assignment of flights to resources is denoted by the binary variables $r_f$ and $c_f$.
Here, $r_f$ denotes if flight $f$ is performed with tail $r$, and $s_f$ denotes if flight $f$ is performed by crew $s$.
Next to flights, the aircraft also have maintenance tasks and reserve duties assigned to them.
Similarly, crew members have scheduled training and leave.
In these periods the crew is not available, and hence we cannot use them for flight activities.
Next to that, we also do not consider the standby (or reserve) duties of crew.
Note that in heavier disruptions, it is safe to assume that no reserve crew is left any more.
This is especially true, when an airline is also struggling with a shortage of flight crew.

\begin{figure*}[ht!]
    \centering
    \input{figure_schedule_example.tex}
    \caption{
    Example schedule consisting of two aircraft and three crew members.
    Here, the horizontal axis represents time, and the flights are drawn as blocks.
    The serrated edge of a block means that the flight departs or lands at an outstation.
    }
    \label{fig:example_schedule}
\end{figure*}
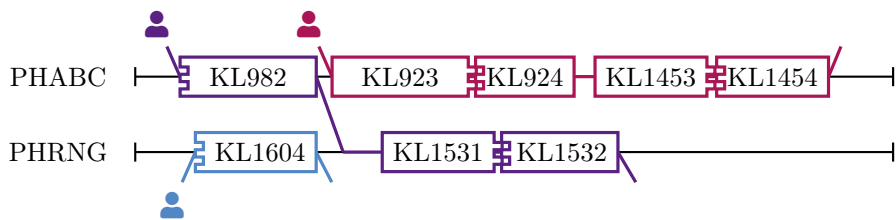

A small example of a schedule is given in \Cref{fig:example_schedule}.
This schedule contains eight flights, two aircraft, and three crew members.
Here, the horizontal axis represents the time of day, and the schedule is shown from the perspective of the fleet assignment.
The flights are represented as blocks, where the straight and serrated edge respectively mark the arrival or departure at the hub-airport or an outstation (non-hub airport).
The horizontal lines represent the assignment of flights for each aircraft, and the coloured lines show the assignment of flights for each crew member.
Note that for the crew assignments, we see that a crew member does not necessarily stay with the aircraft and might need to switch aircraft at some point.
The purple crew member has to switch aircraft between flights KL982 and KL1531.

We look at recovering schedules where many flights experience delays, or where resources are unexpectedly unavailable.
To resolve the disruptions in the schedule, a recovery window is considered in which several levers can be applied.
At the end of the recovery window, the schedule should resume as originally planned.
First, it is possible to make changes to the resource assignment.
That is, change which tail or crew member is assigned to a certain flight.
In practice, this means that we need to swap or move flights between the tails or crew members.
Secondly, we can also make changes to the schedule, that is delaying or cancelling flights.
When delaying a flight, we assume that the flight duration stays the same,
i.e. we do not consider changes in flight speed.
Thus, when a flight $f$ gets delayed by $d$ minutes, the departure and arrival times become $f_d + d$ and $f_a + d$ respectively.
A flight can only be delayed up to a certain maximum $d_\text{max}$, since it is neither possible nor desirable to solve an issue with a huge delay.
Delaying or cancelling a flight has a big impact on the passengers and can consequently be a very expensive lever.

The costs involved can be divided into two categories; namely passenger related costs and operational costs.
Operational costs include fuel cost, maintenance cost, and depreciation of the aircraft.
To represent these, we define $c^{\text{op}}_{f,r}$ to be cost of operating flight $f \in F$ with tail $r \in T$.
Next to operational costs, there may also be passenger related costs for using a certain tail on a certain flight.
This could be for example due to using a smaller plane than in the original plan, meaning that some passengers will need to be reaccommodated as they do not fit on the new aircraft.
We denote this cost with $c^{\text{reassign}}_{f,r}$ for $f \in F$ and $r \in T$.
Note that $c^{\text{reassign}}_{f,r} = 0$, if we perform flight $f$ with the same aircraft type as planned originally.
When a flight is delayed, or cancelled, we incur costs related to passenger claims, rebookings, and, when needed, arranging suitable hotel accommodations.
This is further complicated by passengers on multi-leg journeys.
When one of their flights gets delayed, they might miss their connecting flight, which increases the reaccommodation cost.
These passenger related costs are also called the non-performance cost.

For the non-performance cost, we introduce a binary variable $y_f$ denoting if the flight $f \in F$ is cancelled, and a variable $d_f \in [0, d_\text{max}]$ denoting the delay of flight $f$.
Given a set of passenger connections $P$, that are realized by the undisrupted schedule, we introduce a binary variable $z_p$ to denote if a connection $p \in P$ is broken.
For the cost, there is a function $c_{f}^{\text{delay}} : \mathbb{R} \to \mathbb{R}$ to calculate the delay cost for flight $f$ given the delay in minutes.
We consider this function to be piecewise-constant, making it easy to deal with.
Additionally, the function is increasing.
Lastly, we have parameters $c_{f}^{\text{cancel}}$ and $c_{p}^{\text{mc}}$ to denote the cost of cancelling flight $f$ and breaking passenger connection $p$ respectively.

Summarizing, the total cost is defined as
\begin{equation}\label{eq:objective}
    \sum_{f \in F}\sum_{r \in T} r_f\left(c_{f,r}^{\text{op}} + c_{f,r}^{\text{reassign}}\right)
    + \sum_{f \in F} y_{f} c_{f}^{\text{cancel}}
    + \sum_{f \in F} c_{f}^{\text{delay}}(d_{f})
    + \sum_{p \in P} z_{p} c_{p}^{\text{mc}}
.\end{equation}
Observe that, although counter-intuitive, it is possible that an increase in delay can reduce cost, since it may increase the number of passenger who make a connection.
Also note that in this objective, multiple objectives are combined.
The use of such a single cost objective (or single currency) is similar to current practice within some airlines.

Lastly, a flight assignment is feasible if the following constraints hold:
\begin{enumerate}[(C1)]
    \item\label{constraint:no_overlap}
    \textbf{No Overlap}; There should be enough time between consecutive flights on the same resource.
    Thus, if flight $f$ and $f'$ use the same resource, there needs to be a time $T_{\text{min}}$ between $f$ and $f'$.
    Note that for tails, this  $T_{\text{min}}$ depends on the aircraft type, as bigger aircraft need more time for their ground processes.
    For crew, this minimum connection time depends on the tail assignment for their consecutive duties, as they need less connection time if they do not have to switch aircraft.
    
    \item\label{constraint:rest_and_duty_time}
    \textbf{Rest and Duty Regulations}; The schedule should adhere to the rest and duty times regulations for flight crew.
    We follow European regulations as prescribed by Reg. (EU) No 83/2014.
    In simple terms, a duty consists of a series of flights and is followed by a rest period.
    Duties have a maximum duration based on the number of flights in them.
    The rest period has a minimum length based on the duty before it.
    The exact lengths for these duty and rest periods can be found in the regulatory documents.
    It is not allowed to operate a schedule if these requirements are violated, e.g. delaying a flight such that the end-of-duty time of the crew is exceeded is not feasible, unless we assign another crew.

    \item\label{constraint:tail_crew_match}
    \textbf{Crew Certification}; The assigned flight crew should be certified to fly the assigned aircraft.

    \item\label{constraint:aircraft_flight_match}
    \textbf{Aircraft Authorization}; Aircraft should have the authorization to fly to a certain airport.

    \item\label{constraint:maintenance}
    \textbf{Fixed Maintenance}; The given maintenance tasks and assignments are fixed and cannot be changed.
    Let $M$ be the set of these tasks, then each $m \in M$ is given by a tail $m_t \in T$ on which the maintenance needs to be performed, a start time $m_s$, an end time $m_e$, and an airport where this maintenance is performed $m_l$.
    Then, we have to ensure that $m_t$ is at $m_l$ when the maintenance starts at $m_s$ and cannot be used until $m_e$.

    \item\label{constraint:assign_or_cancel}
    \textbf{Complete Assignment}; A flight $f \in F$ either needs all required resources to be assigned, or the flight is cancelled.
    
    \item\label{constraint:teleportation}
    \textbf{Location}; Lastly, for every resource, the end location for a task in the schedule should be equal to the start location of the next task.
    This constraint implies that we cannot plan any deadheading of aircraft or moving crew as passengers.
    For the last task performed within the recovery window, the end location needs to match the end location in the initial assignment.
    This ensures that it is possible to resume the original schedule after the recovery window is over.
\end{enumerate}

%% file: figure_schedule_example.tex
\begin{tikzpicture}
    \node[] (a_label) at (-1, 0) {PHABC};
    \draw[thick, |-|] (0,0) -- (10, 0);

    \node[flight, inbound,
        minimum width = 1.75cm,
        very thick,
        draw=uuPurple,
    ] (KL982) at (1.5, 0) {KL982};
    
    \node[flight, outbound,
        minimum width = 1.75cm,
        very thick, draw=uuBordeauxRed,
    ] (KL923) at (3.5, 0) {KL923};
    \node[flight, inbound, very thick, draw=uuBordeauxRed]
        (KL924) at (5.14, 0) {KL924};

    \node[flight, outbound, very thick, draw=uuBordeauxRed]
        (KL1453) at (6.8, 0) {KL1453};
    \node[flight, inbound, very thick, draw=uuBordeauxRed] 
        (KL1454) at (8.4, 0) {KL1454};

    \node[] (b_label) at (-1, -1) {PHRNG};
    \draw[thick, |-|] (0, -1) -- (10, -1);

    \node[flight, inbound, very thick, draw=uuBlue]
        (KL1604) at (1.6, -1) {\ KL1604\ };

    \node[flight, outbound, very thick, draw=uuPurple]
        (KL1531) at (4.0, -1) {KL1531};
    \node[flight, inbound, very thick, draw=uuPurple] 
        (KL1532) at (5.6, -1) {\,KL1532};

    \node[color=uuPurple] (crew-1) at (0.3, 0.7) {\faUser};

    \draw[very thick, color=uuPurple] (crew-1)--(KL982.west);
    \draw[very thick, color=uuPurple] (KL982.east)--(2.75,-1);
    \draw[very thick, color=uuPurple] (2.75,-1)--(KL1531.west);
    \draw[very thick, color=uuPurple] (KL1531.east)--(KL1532.west);
    \draw[very thick, color=uuPurple] (KL1532.east)--(6.6, -1.4);

    \node[color=uuBordeauxRed] (crew-2) at (2.3, 0.7) {\faUser};

    \draw[very thick, color=uuBordeauxRed] (crew-2)--(KL923.west);
    \draw[very thick, color=uuBordeauxRed] (KL923.east)--(KL924.west);
    \draw[very thick, color=uuBordeauxRed] (KL924.east)--(KL1453.west);
    \draw[very thick, color=uuBordeauxRed] (KL1453.east)--(KL1454.west);
    \draw[very thick, color=uuBordeauxRed] (KL1454.east)--(9.3,0.4);

    \node[color=uuBlue] (crew-3) at (0.5,-1.7) {\faUser};

    \draw[very thick, color=uuBlue] (crew-3)--(KL1604.west);
    \draw[very thick, color=uuBlue] (KL1604.east)--(2.6,-1.4);

\end{tikzpicture}

%% file: 4_local_search.tex
To solve the problem, we develop a local search algorithm.
Such an algorithm can be tuned to return solutions in a limited time, allowing for the application during the day of operations, where a limited amount of computation time can be afforded.
In order to prevent getting stuck in local optima, we make use of simulated annealing, see \citet{Kirkpatrick1983} and \citet{Cerny1985} for a description of this algorithm.
That is, if the current solution costs $c$, then a new neighbour with lower cost $c'$ is always accepted.
If it has a higher cost, it is accepted with probability $\exp(\frac{c - c'}{T})$, where $T$ is the current temperature that decreases over time. 
Thus, it is possible to accept a `worse' solution, depending on the cost difference and on the current temperature $T$.
As the temperature decreases over time, the probability of accepting worse solutions will decrease as well.
Note that the temperature models the analogy with annealing, i.e. slowly cooling, in physics.

Since we are working with dense schedules, there is not much room for reassigning flights without creating infeasibilities.
For this reason, we increase the search space and relax some hard constraints.
These are the no overlap and the rest and duty times constraints, i.e. Constraints (C\ref{constraint:no_overlap}) and (C\ref{constraint:rest_and_duty_time}) in \Cref{sec:problem_description} respectively.
A violation of these constraints results in a penalty that is added to the cost function.
This penalty is linear in the total time of the violation, i.e. the duration of an overlap for Constraint (C\ref{constraint:no_overlap}), or, for Constraint (C\ref{constraint:rest_and_duty_time}) the length of an extension of duty time or a reduction of rest time beyond their respective limits.
Furthermore, in order to avoid these violations in the final solution as much as possible, the penalty coefficient for such violations is increased during the runtime of the algorithm.

In our model, a disruption consisting of a number of flight delays essentially results in a schedule with lots of violations of the overlap constraint (C\ref{constraint:no_overlap}).
To model a disruption where a resource $r$ is unavailable from time $s$ to $e$, we include a task that starts at $s$, ends at $e$ and is assigned to $r$.
This task cannot be cancelled, or moved, essentially blocking any flight assignments during the period $[s,e]$.
Any flight assigned to $r$ that is in the interval $[s,e]$ will create an overlap with this newly created task, which is expected to be removed by the local search.

In order to reduce the runtime of the algorithm, we steer our local search into only trying to resolve `issues', which are undesirable features currently present in the schedule.
These issues are violations of the relaxed constraints.
Furthermore, missed passenger connections and delayed or cancelled flights are considered as issues as well.
We steer the local search by generating neighbours that affect one of the issues in the current schedule.

An overview of our local search algorithm is given in \Cref{alg:directed_local_search}.
Here, lines \ref{line:get_random_issue}-\ref{line:generate_neighbour} pick a random issue and generate a neighbour to resolve this issue $i$.
This is done by first selecting a random issue and then selecting a flight $f$ to be changed to resolve this issue.
In case of overlapping flights, we always pick the latest of the two.
However, in case of duty time violations, we can pick the flight at the start of the duty or the flight at the end, so we pick at random.
Note that the flight at the start of a duty can only be selected when the duty has not started yet.
Then, the neighbour types that can resolve $i$ by making changes to $f$ are determined, and a neighbour of one of these types is generated.
The neighbour operators are further explained in \Cref{sec:local_search_neighbours}.
Then in \Cref{sec:local_search_neighbour_generation} we explain how to select a neighbour for a given issue and flight.

\begin{algorithm*}
    \caption{Pseudocode for our directed local search algorithm. In here, the variable $c_s$ denotes the cost of solution $s$.}
    \label{alg:directed_local_search}
    \begin{algorithmic}[1]
        \State $s \gets$ Initial schedule with disruptions
        \While{stopping condition is not met}
            \LComment{First select a random issue in the schedule and generate a random neighbour that resolves it.}
            \State $i \gets$ Random issue in the current schedule $s$ 
            \label{line:get_random_issue}
            \State $f \gets$ Flight involved in issue $i$ that can resolve $i$ if changed
            \State Determine neighbour types that can resolve $i$ by changing $f$
            \State Generate a neighbour $n$, of one of these types, that resolves issue $i$
            \label{line:generate_neighbour}
            \State {\color{white}.} % Phantom text so this has proper line height
            \State $s_\text{new} \gets$ Schedule $s$ after applying neighbour $n$
            \LComment{Select the new schedule if it has lower cost or with the probability given by the Boltzmann equation.}
            \If{$c_{s_\text{new}} < c_{s}$}
                \State $s \gets s_\text{new}$
            \ElsIf{\Call{Random}{0,1} $< \exp\left(\tfrac{c_{s} - c_{s_\text{new}}}{T}\right)$}
                \State $s \gets s_\text{new}$
            \EndIf
        \EndWhile
    \end{algorithmic}
\end{algorithm*}

\subsection{Neighbour operators}\label{sec:local_search_neighbours}
The neighbours are derived from the levers we have for resolving a disruption.
We first discuss the neighbours related to changes in the resource assignment, after which we discuss the neighbours for delaying or cancelling flights.

\subsubsection{Changing the Resource Assignment}\label{sec:local_search_neighbours_swaps}
When changing the resource assignment, we either move flights from one resource to another, or swap flights between resources.
Since the connection time of crew members depends on whether they need to switch aircraft, these neighbours do not only affect the assignment, but also the minimum time between tasks.
We have 4 different neighbours to make these assignment changes:
\begin{enumerate}
    \item \textbf{SwapTail}; This neighbour swaps two series of flights between two tails.
    \item \textbf{SwapCrew}; Similar to SwapTail, but swapping between crew.
    \item \textbf{SwapBoth}; Swaps two series of flights between two tail-crew pairs. Thus, in this neighbour, the tail and crew stay together.
    \item \textbf{MoveTail}; Moves a series of flights from one tail to another.
\end{enumerate}

In these neighbours, we first need to select the flights to be swapped or moved.
As mentioned before, before a neighbour is generated, a specific issue $i$ and flight $f$ are selected.
However, we cannot simply swap or move a single flight, as this would most likely create violations of Constraint (C\ref{constraint:teleportation}), the location constraint.
To ensure that this constraint is not violated, we consider rotations of flights.
A rotation is a series of flights starting at a hub airport $h$ and ending at the same airport, without going to $h$ in between.

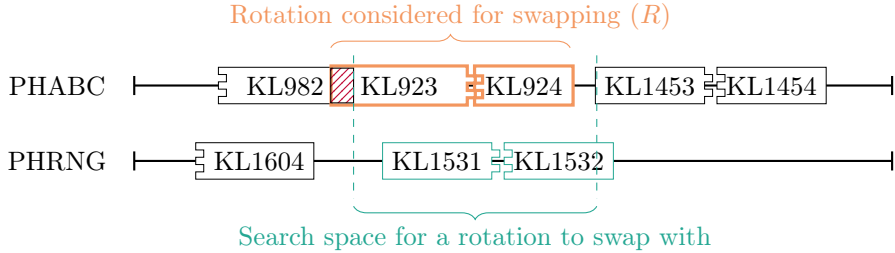
\begin{figure*}[ht!]
    \centering
    \input{figure_swap_neighbour.tex}
    \caption{
    Illustration of the swap neighbour.
    Here, there is an overlap between the flights KL982 and KL923, which we try to fix.
    The rotation considered for moving ($R$) is highlighted orange (KL923 and KL924), and we try to swap this with the green highlighted flights (KL1531 and KL1532).
    }
    \label{fig:swap_neighbour_illustration}
\end{figure*}

First, we look at the SwapTail neighbour.
This neighbour considers the tail $t$ assigned to flight $f$, and selects a random tail $t'$ of the same type to swap with.
Then, we consider the rotation $R$ containing flight $f$.
Note that $R$ may contain flights that cannot be reassigned because they are already underway.
Furthermore, it is also possible for $R$ to contain a maintenance task.
In these cases, the neighbour will fail to generate.
Let $f^{\text{prev}}_a$ denote the arrival time of the first flight before $R$ of $t$, and $f^{\text{next}}_d$ denotes the departure time of the first flight after $R$ on $t$.
To determine the flights on $t'$ to swap with, we search close in time to $R$.
That is, we search for a rotation $R'$ that departs in the interval $[f^{\text{prev}}_a, f^{\text{next}}_d]$.
As with $R$, the rotation $R'$ may not contain flights that cannot be reassigned or maintenance tasks.
If no such rotation is found, the neighbour fails to generate.
The neighbour is visualized in \Cref{fig:swap_neighbour_illustration}.
In this image, we see the flight assignment of two aircraft.
Note that this is the same assignment as the example in \Cref{sec:problem_description}, but the crew assignment is not visualized here.
Furthermore, flight KL982 incurred a delay, which resulted in an overlap between the flights KL982 and KL923.
This is the issue that we try to fix.
The serrated edge on a flight block illustrates that the arrival or departure on a non-hub airport.
Note that in this particular example, the selected swap will create another overlap on flights KL1532 and KL1453.
For this reason, we include the option to extend the selection of rotations to swap.
With a given probability ($60\%$ in our case), we try to extend the selected rotations to swap (currently $R$ and $R'$) by adding their respective next rotations until both aircraft are at the same location at the same time.
If successful, this variant will not create any additional overlaps.
In this particular example, that would mean that not only $R$ is selected, but also the flights KL1453 and KL1454.

The SwapCrew neighbour is very similar to the SwapTail neighbour.
It selects flights in the same way, and the crew members selected for the swap need to have the same qualifications as the crew assigned to $f$.
Furthermore, the cost calculations on this neighbour are more elaborate because in this case we need to check and recalculate the duty and rest times of the impacted crew.

In the last swapping neighbour, SwapBoth, we keep crew and tail together.
Selecting flights works similarly to the previous two, however, this time we need to check that all flights in the selected series are performed by the same tail-crew pair.
When selecting a new tail, it is now possible to select a tail with a different type.
This is possible in this variant, since the crew stays with their assigned aircraft, hence, it is guaranteed that the crew is licensed for the aircraft they operate.
The possibility of doing these type swaps adds extra flexibility when recovering, and was rarely included in previous literature.

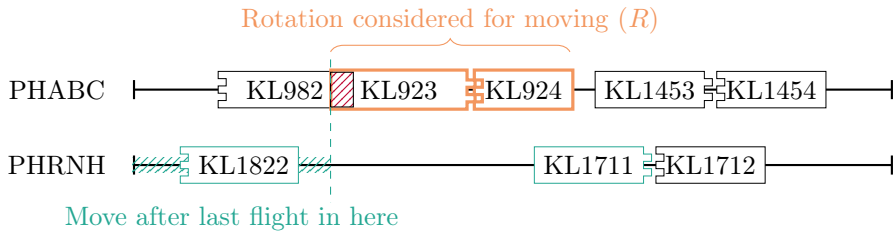
\begin{figure*}[ht!]
    \centering
    \input{figure_move_neighbour.tex}
    \caption{
    Illustration of the move neighbour.
    Here, there is an overlap between the flights KL982 and KL923, which we try to fix.
    The rotation considered for moving ($R$) is highlighted orange (KL923 and KL924), and we try to move this in between the green highlighted flights (KL1822 and KL1711).
    }
    \label{fig:move_neighbour_illustration}
\end{figure*}

Lastly, there is the MoveTail neighbour.
Here, we select the rotation $R$ as before, but instead of swapping, this rotation is moved to a different tail.
As in the SwapTail neighbour, this new tail is a random tail $t'$ of the same type as the one assigned to $f$.
On this new tail, we need to determine where to move $R$.
For this, let $s$ denote the departure time of $R$.
Note that this is different from the previously defined $f^{\text{prev}}_a$.
Then $R$ is moved after the last flight of tail $t'$ that arrives at $h$ before time $s$.
This may create an overlap if the tail already has flights assigned around the same time as $R$.
However, if this is acceptable or not is for the simulated annealing to determine.
An example of this neighbour is illustrated in \Cref{fig:move_neighbour_illustration}.
Here, we consider the same issue as in \Cref{fig:swap_neighbour_illustration}, but a different aircraft is selected to move these flights to.

Note that with these neighbours, our algorithm currently only considers a single hub airport $h$, however, with some minor changes the neighbours should also work with multi-hub networks.
For this, one needs to make sure that when moving flights, the start and end location needs to be the same airport.
When swapping flights, this is not the case, as we can swap series of flights if both series have the same start and end airport, respectively.

\subsubsection{Changing the Schedule}\label{sec:local_search_neigbhours_schedule}
Apart from changing the resource assignment, it is also possible to make changes in the schedule.
For this, we have neighbours that delay or cancel flights.
Since the local search performs numerous iterations, there could be situations where it becomes possible to improve the solution by removing the delay from a flight that was delayed earlier in the search process.
The same is true for cancelled flights.
Thus, there are also neighbours that try to undo delays and cancellations.

For the delay neighbour, we calculate exactly how much flight $f$ needs to be delayed in order to resolve the selected issue $i$.
Next to delaying $f$, the neighbour might also delay future flights that are in the same rotation as $f$.
Issues, such as overlapping tasks, occurring at an outstation have limited recovery options.
It is not possible to fix these using the previously mentioned reassignment neighbours, as there are no resources to swap with.
Thus, if an outbound flight is delayed such that it overlaps with the inbound flight, the neighbour automatically delays this inbound flight as well.
Note here that there might be some buffer between these two flights, so the delays of these flights are not necessarily the same.

The next schedule change is to cancel a flight.
While this undesirable for the passengers, it frees up resources such that they can be used elsewhere.
Note that cancelling flights is quite an expensive action, which may lead to difficulties accepting this neighbour.
However, when cancelling a flight also resolves a soft constraint violation, the penalty reduction can compensate the cancellation cost.
The cancel neighbour cannot simply cancel a single flight, as this would violate Constraint (C\ref{constraint:teleportation}) and leave our resources stranded.
For this reason, only full rotations can be cancelled.
Thus, this neighbour proposes to cancel the rotation containing flight $f$.
However, this cancellation is not always possible.
If a crew change occurs during the rotation and the rotation is cancelled, the new crew cannot come back to $h$ and resume their duty.
Thus, in this case, the neighbour fails to generate.

\subsection{Neighbour Generation}\label{sec:local_search_neighbour_generation}
As mentioned before, we perform directed search in the sense that we steer the local search to try to resolve `issues' currently present in the schedule.
This is done by randomly selecting an issue and generating one or more neighbours that could resolve it.
Note that resolving one issue, might create other issues.
After selecting an issue as described in \Cref{alg:directed_local_search}, we determine which neighbour types can resolve this issue.
This is done by following a lookup-table.
In case of a delayed flight, only the removal of this delay can be tried, and similar for cancelled flights.
In the other cases, the delay and cancel neighbours are always possible neighbour types.
To determine which reassignment neighbours to generate, we use \Cref{tab:lookup_table_swap_generation}.
We denote these possible neighbour types with the set $N$.
With these types determined, a specific neighbour is generated.
For this last step, we try different strategies.

\begin{table*}[ht!]
    \centering
    \caption{Lookup table to determine which reassignment neighbours to generate for a given issue type.}
    \label{tab:lookup_table_swap_generation}
    \begin{tabular}{l cccc}
        \toprule
        & SwapTail & SwapCrew & SwapBoth & MoveTail \\
        Issue & & & & \\
        \midrule
        Overlap Tail & \checkmark & & \checkmark & \checkmark \\
        Overlap Crew & & \checkmark & \checkmark & \\
        Overlap Tail \& Crew & \checkmark & \checkmark & \checkmark & \checkmark \\
        Exceeded Duty & & \checkmark & \checkmark & \\
        Not Enough Rest & & \checkmark & \checkmark & \\
        \bottomrule
    \end{tabular}
\end{table*}

The first strategy closely follows the more `random' approach that is classical for local search algorithms.
Thus, we select a random $n \in N$, and generate a neighbour of this type.
Recall that the swap and move neighbourhoods randomly select a tail or crew to respectively swap with or move to.
The other neighbours are deterministic given an issue.
For this random selection, we set a weight for each neighbour type.
This approach will be denoted as ``RandomNeighbour''.

In the second approach, we generate multiple candidate neighbours and select a neighbour from this set.
Thus, for each $n \in N$, we generate $25$ neighbours of type $n$.
Note, however, that the delay and cancel neighbour are deterministic, thus for those only one candidate is generated.
The set containing these generated neighbours is called $C$.
Then, we need to select one $c \in C$ to be applied.
We test multiple methods to make such a selection.
The first selection method, denoted as ``CandidatesBest'', is to simply pick the neighbour that leads to the greatest improvement.
This selection method, however, might be too greedy.
Thus, we also use a method that tries to balance this greedy approach and the random approach discussed before.
This method is denoted as ``CandidatesSoftMin'', where, instead of selecting the best neighbour, we select a random neighbour according to a probability as given by a so-called SoftMin function.
If we have a set of candidate neighbours $C$, and a function $\delta : C \to \mathbb{R}$ to calculate the cost difference between a solution and a neighbour in $C$.
Then, the probability of selecting a neighbour $c \in C$ is given by
\begin{equation}\label{eq:softmin_neighbour}
    \mathds{P}(\text{choosing neighbour } c)
    = \frac{
        e^{-\delta(c)}
    }{
        \sum_{c' \in C} e^{-\delta(c')}
    }
.\end{equation}
This method scales the delta-scores in such a way that better neighbours have a higher probability of being selected.

A third strategy we considered was a variable neighbourhood approach.
Here we first consider issues affecting the aircraft for a given number of iterations, and perform a number of iterations where we focus on issues affecting crew, and so on.
However, early testing showed that this did not work as well as the other approaches, thus we excluded it in our final tests.

%% file: figure_swap_neighbour.tex
\begin{tikzpicture}
    \node[] (a_label) at (-1, 0) {PHABC};
    \draw[thick, |-|] (0,0) -- (10, 0);

    \node[flight, inbound,
        minimum width = 1.75cm
    ] (KL982) at (2.0, 0) {KL982};
    
    \node[flight, outbound, minimum width = 1.75cm, selected]
        (KL923) at (3.5, 0) {KL923};
    \node[flight, inbound, selected]
        (KL924) at (5.14, 0) {KL924};

    \node[flight, outbound] (KL1453) at (6.8, 0) {KL1453};
    \node[flight, inbound]  (KL1454) at (8.4, 0) {KL1454};

    \fill[pattern color=uuRed, pattern=north east lines, draw=black] (2.6, 0.23) rectangle (2.9, -0.23) {};

    \node[] (b_label) at (-1, -1) {PHRNG};
    \draw[thick, |-|] (0, -1) -- (10, -1);

    \node[flight, inbound] (KL1604) at (1.6, -1) {\ KL1604\ };

    \node[flight, outbound, draw=uuGreen] (KL1531) at (4.0, -1) {KL1531};
    \node[flight, inbound, draw=uuGreen]  (KL1532) at (5.6, -1) {KL1532};

    \draw[decorate, uuOrange, decoration={brace,amplitude=5pt}]
        (2.6, 0.5) -- (5.75, 0.5) node[midway,yshift=1.2em]{Rotation considered for swapping ($R$)};

    \draw[uuGreen, dashed] (2.9,0.4) -- (2.9,-1.5);
    \draw[uuGreen, dashed] (6.1,0.4) -- (6.1,-1.5);
    \draw[decorate, uuGreen, decoration={brace,amplitude=5pt,mirror}]
        (2.9, -1.6) -- (6.1, -1.6) node[midway,yshift=-1.2em]{Search space for a rotation to swap with};
\end{tikzpicture}

%% file: figure_move_neighbour.tex
\begin{tikzpicture}
    \node[] (a_label) at (-1, 0) {PHABC};
    \draw[thick, |-|] (0,0) -- (10, 0);

    \node[flight, inbound,
        minimum width = 1.75cm
    ] (KL982) at (2.0, 0) {KL982};
    
    \node[flight, outbound, minimum width = 1.75cm, selected]
        (KL923) at (3.5, 0) {KL923};
    \node[flight, inbound, selected]
        (KL924) at (5.14, 0) {KL924};

    \node[flight, outbound] (KL1453) at (6.8, 0) {KL1453};
    \node[flight, inbound]  (KL1454) at (8.4, 0) {KL1454};

    \fill[pattern color=uuRed, pattern=north east lines, draw=black] (2.6, 0.23) rectangle (2.9, -0.23) {};

    \begin{scope}[on background layer]
        \node[] (b_label) at (-1, -1) {PHRNH};
        \draw[thick, |-|] (0, -1) -- (10, -1);
    \end{scope}

    \node[flight, inbound, draw=uuGreen] (KL1822) at (1.4, -1) {\ KL1822\ };

    \node[flight, outbound, draw=uuGreen] (KL1711) at (6.0, -1) {KL1711};
    \node[flight, inbound]                (KL1712) at (7.6, -1) {KL1712};

    \begin{scope}[on background layer]
        \fill[pattern color=uuGreen, pattern=north east lines, on background layer] (0, -0.9) rectangle (2.6, -1.1)
            node[midway, yshift=-2.0em, color=uuGreen]{Move after last flight in here};        
    \end{scope}
    \draw[uuGreen, dashed] (2.6,0.4) -- (2.6,-1.5);

    \draw[decorate, uuOrange, decoration={brace,amplitude=5pt}]
        (2.6, 0.5) -- (5.75, 0.5) node[midway,yshift=1.2em]{Rotation considered for moving ($R$)};

\end{tikzpicture}

%% file: 5_results.tex
In this section, we present our computational study to analyse the performance of our local search algorithm.
We consider both the runtime and the solution quality and compare the different neighbour generation methods.
Furthermore, we investigate the effect of different cost choices in the model.
These are changes to the priority of passengers missing their connection with respect to delays, essentially changing the delay cost curve.
A penalty for changing resource assignments is also considered.
Lastly, we compare our integrated method to sequential solution approaches.
We first discuss the setup of the experiments, after which present the results.

\subsection{Setup}\label{sec:results_setup}
We test the algorithm with data from KLM Royal Dutch Airlines, making use of one of their recent flight networks as well as their cost calculations.
These cost estimations are disturbed to avoid revealing confidential data.
Our test instance is the European flight network of KLM, consisting of roughly \num{500} flights per day.
This network has a single hub airport and is served by \num{6} different aircraft types.
The recovery window is $72$ hours, meaning that from now, there is a $72$ hour window in which we can apply changes to recover from the current disruption.
This is done in order to get a good overview of the propagation effects and their cost.
As mentioned in \Cref{sec:problem_description}, we analyse two types of disruptions, namely flight delays (FD) and unavailable flight resources (UR).
We consider the occurrence of disruptions at different times of the day: 8:00, 12:00, 15:00, and 18:00.
The time of day affects both the number of options and the kind of options we have for resolving disruptions.
For example, later in the day the crew has, in general, less duty time remaining, hence delay options might be limited.
On the other hand, if we look at unavailable resources, it might be a lot `easier' if they occur later in the day, as they will not affect that many flights.

For our tests, we generate $25$ instances for each disruption type and time pair.
To generate an instance of the flight delay disruption, each flight that is currently underway has a $25\%$ probability to get a random delay.
These flights have an arrival delay of 30 minutes on average.
This results in 25 instances of a disrupted schedule that, on average, contain about \num{29} overlapping flight assignments and a total overlap of about $469$ minutes.
Note this represents a heavily disrupted day, which is not a daily occurrence.
However, these are the kind of disruptions where these algorithms are needed the most.
For the resource unavailability, we select one or more flights that just landed or are close to departing.
Then, either the crew, aircraft, or both become unavailable and we block the disrupted resource for the next $6$ to $12$ hours, meaning that this resource is unavailable for those $6$ to $12$ hours.
On average, this results in about \num{4} flights being disrupted.

When we compare the quality of the solutions, we look at a few performance indicators.
For the cost of a solution, we compare both the total cost and the Non-Performance Cost (NPC).
The NPC denotes the sum of the delay, cancellation, and missed passenger costs, and essentially tells something about the impact on the passengers.
The total cost of a solution is the NPC together with the operational costs, thus the NPC is the sum of the last three terms in \Cref{eq:objective}.
Furthermore, we compare the number of cancelled flights, number of delayed flights, and number of missed passenger (pax) connections.
Here, we specifically highlight the number of flights that have a delay of more than 15 minutes.
Lastly, we also compare the maximum and average delays in the schedule, where we calculate the average delay over only the delayed flights.

We run our local search in a parallel multi-start setup.
That is, we start multiple parallel searches to see which one of them comes up with the best solution.
The only difference between these searches is the seed used in the random generator.
As a stopping criterion for the local search, we look at the number of iterations since the last accepted neighbour.
When this exceeds $10\ 000$, we either restart or stop.
In total, the search is restarted twice, where each time the temperature is increased.
To account for the stochasticity in the search procedure, we use the 5 best solutions in our comparisons.
Hence, to calculate a KPI value for a single disruption, we take the average over the 5 best solutions returned by the simulated annealing.

\subsection{Comparing Neighbour Generation Strategies}\label{sec:results_neigbhour_generation_strategies}
First, we compare the neighbour generation methods described in \Cref{sec:local_search_neighbour_generation}.
For this, we present the average cost change for using a candidate method compared to the RandomNeighbour method.
This is shown in \Cref{tab:search_strategy_compare_cost}.
For each combination of disruption type and time of day, we take the average over the 25 different disruption instances.
Note that our model looks 3 days ahead and considers operational cost of the schedule for the complete European flight network of KLM.
This means that the total cost of a solution is quite large, and that small percentage improvements lead to big cost reductions.
From this table, we see that in the flight delay scenario, the candidate generation methods generally perform better, although this advantage becomes smaller later in the day.
Still, the use of one of our candidate generation methods should lead to better results in this scenario.
When resources are unavailable, there is not a clear difference between the different methods.

\begin{table*}[ht!]
\centering
\caption{Average percentage cost change of each neighbour generation strategy compared to the RandomNeighbour method.}
\label{tab:search_strategy_compare_cost}
\begin{tabular}{llrrrr}
\toprule
 &  & 08:00 & 12:00 & 15:00 & 18:00 \\
Disruption & Search strategy &  &  &  &  \\
\midrule
\multirow[c]{2}{*}{FD} & CandidatesBest & \num{-0.36}\% & \num{-0.25}\% & \num{-0.17}\% & \num{-0.04}\% \\
 & CandidatesSoftMin & \num{-0.33}\% & \num{-0.20}\% & \num{-0.15}\% & \num{-0.04}\% \\
\cmidrule[0.01em]{1-6}
\multirow[c]{2}{*}{UR} & CandidatesBest & \num{0.00}\% & \num{-0.01}\% & \num{-0.03}\% & \num{-0.03}\% \\
 & CandidatesSoftMin & \num{-0.05}\% & \num{-0.07}\% & \num{-0.03}\% & \num{-0.02}\% \\
\bottomrule
\end{tabular}
\end{table*}

In the detailed KPI comparison in \Cref{tab:search_strategy_compare_kpi}, we see a similar picture.
Here, the candidate methods also show a better performance in case of flight delays.
Compared to the RandomNeighbour method, they show a lower overall cost as is also seen in \Cref{tab:search_strategy_compare_cost}.
These candidate methods also show a reduction in NPC.
However, in case of unavailable resources, the picture becomes different again.
Here, we see increases in the NPC, while the overall cost is lower.
In this scenario, there are fewer delays and cancellations, thus choosing a different flight to delay or cancel can have a relatively big impact on the NPC.
This is due to the cost differences between different flights.
Combined with the randomness of the search procedure and the NPC being only part of the objective, this creates some variation in the NPC.
However, the total cost, which is the objective, does not show this behaviour, showing that the higher NPC does allow for operational cost savings.
From this perspective, the CandidatesBest method is best in case of flight delays, while the CandidatesSoftMin method performs better when resources become unavailable.

\begin{table}[ht!]
\centering
\caption{Comparison of average KPI values for different neighbour generation strategies.}
\label{tab:search_strategy_compare_kpi}
\begin{tabular}{llrrr}
\toprule
 &  & RandomNeighbour & CandidatesBest & CandidatesSoftMin \\
Disruption & KPI &  &  &  \\
\midrule
\multirow[c]{8}{*}{FD} & \# Cancelled flights & 0.12 & 0.12 & 0.15 \\
 & \# Delayed flights & 59.39 & 50.13 & 53.64 \\
 & \# Missed pax connections & 52.08 & 51.31 & 50.88 \\
 & \# $> 15$ min delayed & 17.97 & 14.80 & 15.20 \\
 & Maximum delay (minutes) & 54.38 & 53.51 & 53.85 \\
 & Average delay (minutes) & 14.08 & 13.70 & 13.62 \\
 & NPC change & n.a. & -10.08\% & -8.93\% \\
 & Total cost change & n.a. & -0.20\% & -0.17\% \\
\cmidrule[0.01em]{1-5}
\multirow[c]{8}{*}{UR} & \# Cancelled flights & 0.34 & 0.29 & 0.32 \\
 & \# Delayed flights & 2.52 & 1.71 & 1.58 \\
 & \# Missed pax connections & 2.37 & 2.40 & 2.19 \\
 & \# $> 15$ min delayed & 0.66 & 0.62 & 0.52 \\
 & Maximum delay (minutes) & 38.54 & 30.57 & 28.02 \\
 & Average delay (minutes) & 15.28 & 20.09 & 14.98 \\
 & NPC change & n.a. & +76.61\% & +81.48\% \\
 & Total cost change & n.a. & -0.02\% & -0.04\% \\
\bottomrule
\end{tabular}
\end{table}

However, this comes at a performance cost, as can be seen in \Cref{tab:search_strategy_compare_runtime}.
Here, we compare the average runtime of the different neighbour generation methods.
The RandomNeighbour strategy is a lot faster than the candidate selection methods, running in just a few seconds.
This is as expected, since the candidate selection methods generate more neighbours in a single iteration, thus they perform fewer iterations per second.
We further show this in \Cref{tab:search_strategy_compare_iterations}, where the average number of iterations of the different neighbour generations methods is compared.
Comparing the runtime with the solution quality in \Cref{tab:search_strategy_compare_cost}, we find that the RandomNeighbour method still finds quite a good solution in such a short amount of time.
Thus, when time becomes critical, we can get a good solution almost instantly.
Note that in a multi-start approach, it is possible to use multiple neighbour generation methods in different rounds.
Thus, a definitive choice of neighbour generation method is not strictly necessary.

\begin{table*}[ht!]
\centering
\caption{Average runtime (in seconds) of a single simulated annealing run.}
\label{tab:search_strategy_compare_runtime}
\begin{tabular}{llrrrr}
\toprule
 &  & 08:00 & 12:00 & 15:00 & 18:00 \\
Disruption & Search strategy &  &  &  &  \\
\midrule
\multirow[c]{3}{*}{FD} & CandidatesBest & 7.900 & 24.909 & 18.449 & 19.807 \\
 & CandidatesSoftMin & 7.216 & 25.118 & 19.443 & 20.539 \\
 & RandomNeighbour & 6.929 & 4.053 & 3.768 & 3.753 \\
\cmidrule[0.01em]{1-6}
\multirow[c]{3}{*}{UR} & CandidatesBest & 7.659 & 5.134 & 2.809 & 2.331 \\
 & CandidatesSoftMin & 10.810 & 8.048 & 5.936 & 6.534 \\
 & RandomNeighbour & 2.408 & 1.879 & 1.476 & 1.233 \\
\bottomrule
\end{tabular}
\end{table*}

\begin{table*}[ht!]
\centering
\caption{Average number of iterations (in millions) of a single simulated annealing run.}
\label{tab:search_strategy_compare_iterations}
\begin{tabular}{llrrrr}
\toprule
 &  & 08:00 & 12:00 & 15:00 & 18:00 \\
Disruption & Search strategy &  &  &  &  \\
\midrule
\multirow[c]{3}{*}{FD} & CandidatesBest & 3.03 & 2.15 & 1.91 & 1.94 \\
 & CandidatesSoftMin & 2.41 & 2.20 & 1.96 & 2.56 \\
 & RandomNeighbour & 3.12 & 2.17 & 1.90 & 2.36 \\
\cmidrule[0.01em]{1-6}
\multirow[c]{3}{*}{UR} & CandidatesBest & 0.61 & 0.20 & 0.13 & 0.11 \\
 & CandidatesSoftMin & 0.87 & 0.54 & 0.49 & 0.48 \\
 & RandomNeighbour & 0.92 & 0.48 & 0.41 & 0.32 \\
\bottomrule
\end{tabular}
\end{table*}

\subsection{Comparing Solution Quality}\label{sec:results_solution_quality}
Next, we look at the quality of the solutions found by the local search.
For this, we use the CandidatesBest neighbour generator for FD disruptions and CandidatesSoftMin for UR disruptions.
To assess the solution quality, we compare our local search with a delay propagation algorithm and with different sequential approaches.
The comparison with sequential approaches is to compare and show the effect of optimal tail recovery on the remaining crew problem, thus overall computation time of these algorithms is less relevant.
In current literature not many integrated approaches consider a similar crew model to ours, with many of them being simplified.
Moreover, the few models that are available either have a very large runtime, or do not publish sufficient information on their algorithm and runtime.
Therefore there are no clear state-of-the-art methods to directly compare to.

First compare with a delay propagation algorithm, as this would show us what happens when the minimal actions are taken when resolving a disruption.
This algorithm resolves disruptions by propagating delays or cancelling flights.
It performs the following steps:
\begin{enumerate}
    \item Sort the issues in the schedule from earliest to latest. In here, we only consider the violations of the relaxed constraints (Constraints (C\ref{constraint:no_overlap}) and (C\ref{constraint:rest_and_duty_time})).
    \item Pick the issue that starts the earliest.
    \item Try to resolve the issue by delaying subsequent flights. If this cannot be done, resolve the issue by cancelling the flight rotation.
    \item Repeat these steps until all issues are resolved.
\end{enumerate}
Note that there could be situations where this approach does not return a feasible solution, i.e. it can not resolve all issues present in the schedule.
An example of this situation is when the delay of a flight leads to exceeded duty hours of the crew, but cancelling is also not possible due to a crew change on the outstation.
In our presentation, we only include the cases where the delay propagation returns a feasible solution.
Note, however, that our local search can solve more situations than this simple benchmarking algorithm.
In our results, this algorithm is denoted as $P_{\text{int}}$.

Next to comparing to this delay propagation algorithm, we also make comparisons to sequential solutions.
For this, we first compare to a method that we call `crew-naive'.
The crew-naive method uses the local search to resolve all tail issues, but uses delay propagation to resolve the remaining crew issues.
This algorithm will be denoted as $LS_{t}|P_{c}$, for doing local search on the tail schedule and use a delay propagation solver for the crew.

We also compare to a method that is commonly used in practice, where the tail issues are resolved optimally.
For this, we use the MIP introduced in \Cref{sec:tail_mip}.
The remaining crew issues are solved using the previously described delay propagation method.
We denote this method as $MIP_t|P_c$.
The MIP model has a limited horizon in order to reduce the size, looking only 36 hours ahead.
This reduction in horizon was necessary in order to still get results within a reasonable time-frame.
The MIP is solved by Gurobi 12.0.0, with some slight parameter changes.
The parameter \texttt{MIPFocus} is set to finding feasible solutions, heuristics are increased to $30\%$, and a 30-minute runtime restriction is imposed.
Nevertheless, the optimality gap is at most $0.4\%$, and on average $0.03\%$.
Note that both our sequential algorithms are not well suited to handle scenarios containing only crew issues.
In these cases one would expect the solution to be more or less the same as $P_\text{int}$.

\begin{figure*}[ht!]
    \centering
    \includegraphics[width=0.85\textwidth]{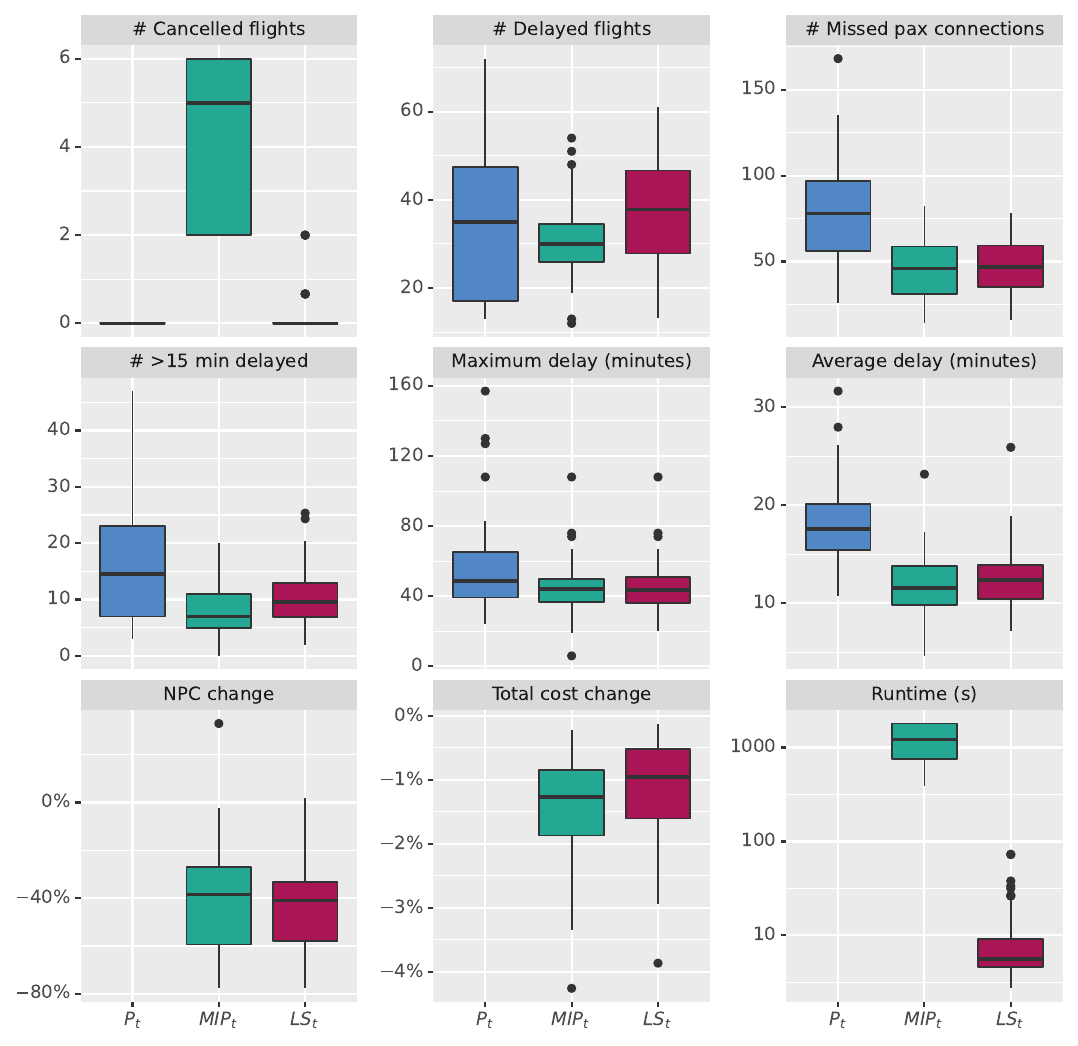}
    \caption{Comparison of KPI values when recovering the tail schedule in case of the FD scenario.}
    \label{fig:tail_kpi_boxplot_fd}
\end{figure*}

\begin{figure*}[ht!]
    \centering
    \includegraphics[width=0.85\textwidth]{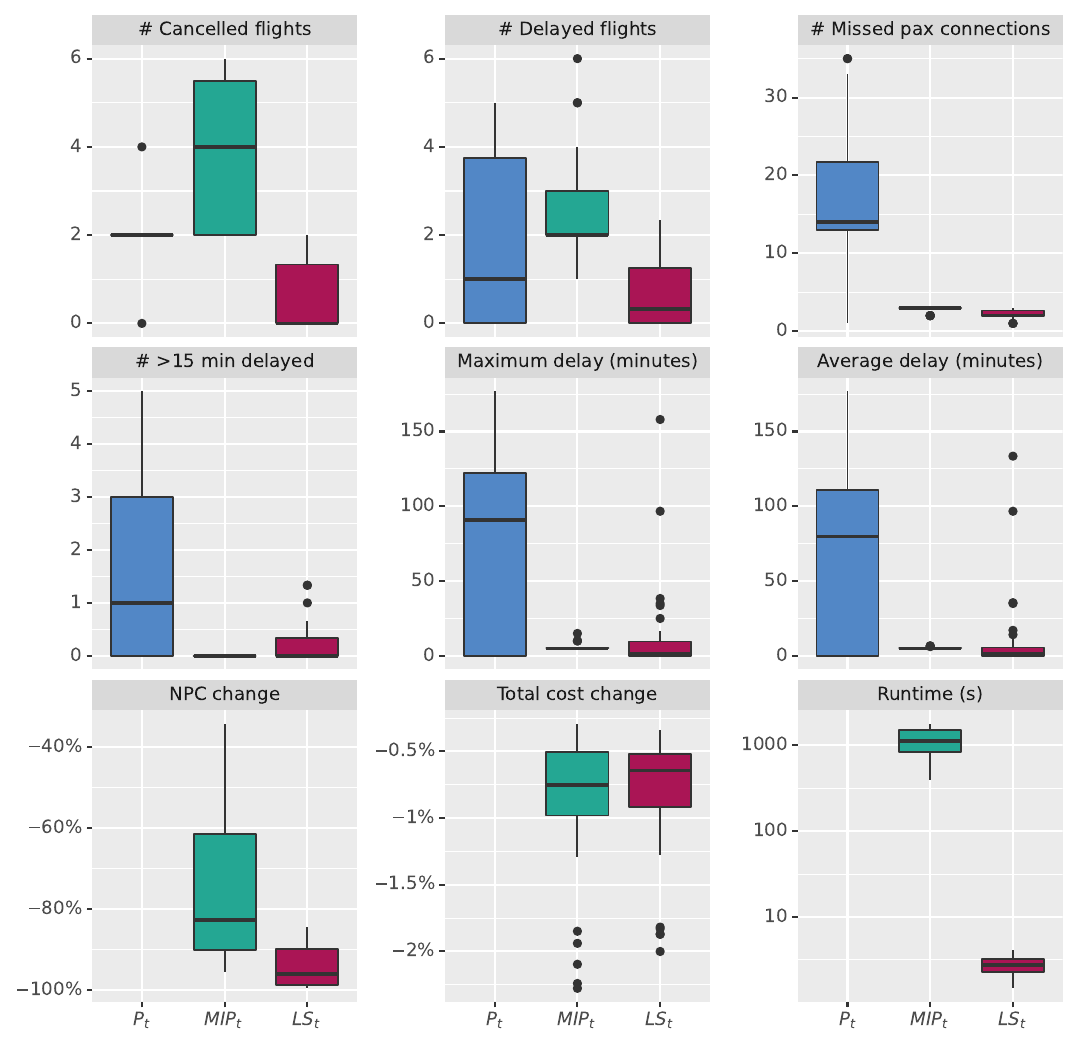}
    \caption{Comparison of KPI values when recovering the tail schedule in case of the UR scenario.}
    \label{fig:tail_kpi_boxplot_ur}
\end{figure*}

To compare different types of algorithms, we first ran a comparison doing only tail recovery.
Here, we compare the solution quality of the MIP ($MIP_t$) to the above-mentioned delay propagation method ($P_t$) and our local search ($LS_t$).
For this, we excluded all scenarios containing only crew issues and evaluate on only the tail assignment.
We show this comparison in \Cref{tab:tail_quality_kpi}.
In this table, we show the average values for each KPI.
We also show how these values are distributed over the different scenarios in \Cref{fig:tail_kpi_boxplot_fd,fig:tail_kpi_boxplot_ur}.
These show the KPI values in case of flight delays and unavailable resources respectively.
In here, we see that $MIP_t$ indeed shows the best performance, as it has the lowest total cost.
The local search, however, shows close to optimal solutions while taking significantly less computation time to find them.
Moreover, if we look at the NPC, we see that $LS_t$ shows better results than the MIP, especially in the UR scenario.
The MIP compensates for this by having a lower operational cost, resulting in a lower total cost.
This could be done by, for example, choosing more cost-effective reassignments, however, it seems that $MIP_t$ needs to cancel more flights in order to achieve this goal.
Local Search is usually quite hesitant to choose cancellation, i.e. it selects this neighbourhood with a small probability and sees cancellations as issues it wants to resolve.

\begin{table*}[ht!]
\centering
\caption{Comparison of average KPI values when recovering the tail schedule.}
\label{tab:tail_quality_kpi}
\begin{tabular}{llrrr}
\toprule
 &  & $P_t$ & $MIP_t$ & $LS_t$ \\
Disruption & KPI &  &  &  \\
\midrule
\multirow[c]{9}{*}{FD} & \# Cancelled flights & 0.00 & 4.31 & 0.15 \\
 & \# Delayed flights & 34.81 & 30.84 & 37.62 \\
 & \# Missed pax connections & 78.90 & 44.57 & 46.53 \\
 & \# $> 15$ min delayed & 16.69 & 7.78 & 10.06 \\
 & Maximum delay (minutes) & 56.17 & 44.45 & 45.39 \\
 & Average delay (minutes) & 18.34 & 11.78 & 12.61 \\
 & NPC change & n.a. & -40.20\% & -43.71\% \\
 & Total cost change & n.a. & -1.45\% & -1.17\% \\
 & Runtime (s) & n.a. & 1241.77 & 11.58 \\
\cmidrule[0.01em]{1-5}
\multirow[c]{9}{*}{UR} & \# Cancelled flights & 2.00 & 3.71 & 0.57 \\
 & \# Delayed flights & 1.82 & 2.76 & 0.55 \\
 & \# Missed pax connections & 15.26 & 2.82 & 2.13 \\
 & \# $> 15$ min delayed & 1.71 & 0.00 & 0.22 \\
 & Maximum delay (minutes) & 75.50 & 5.74 & 13.65 \\
 & Average delay (minutes) & 70.22 & 5.17 & 11.56 \\
 & NPC change & n.a. & -75.29\% & -94.23\% \\
 & Total cost change & n.a. & -0.91\% & -0.87\% \\
 & Runtime (s) & n.a. & 1122.66 & 2.78 \\
\bottomrule
\end{tabular}
\end{table*}

\begin{figure*}[ht!]
    \centering
    \includegraphics[width=0.85\textwidth]{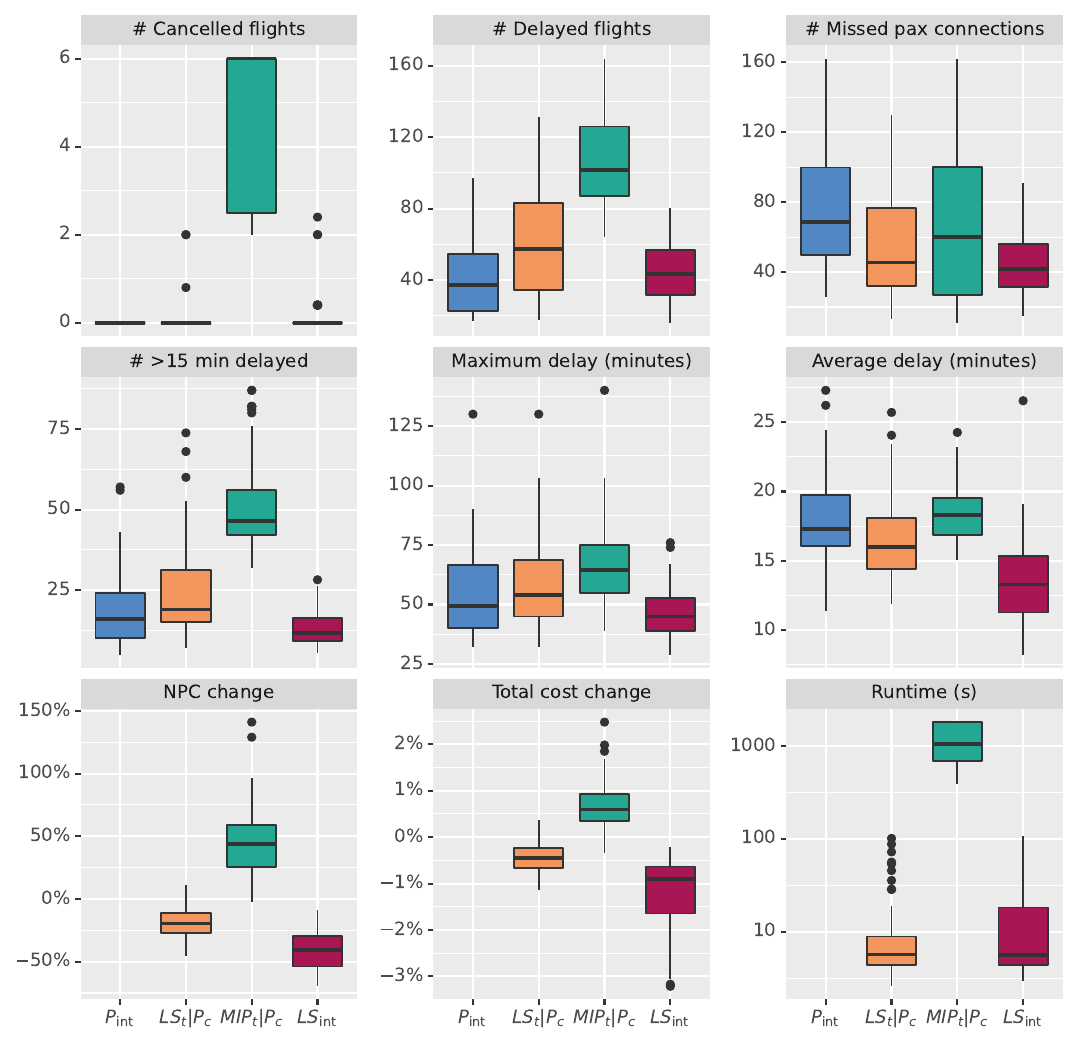}
    \caption{Comparison of KPI values when recovering both tail and crew schedules in case of the FD scenario.}
    \label{fig:kpi_boxplot_fd}
\end{figure*}

\begin{figure*}[ht!]
    \centering
    \includegraphics[width=0.85\textwidth]{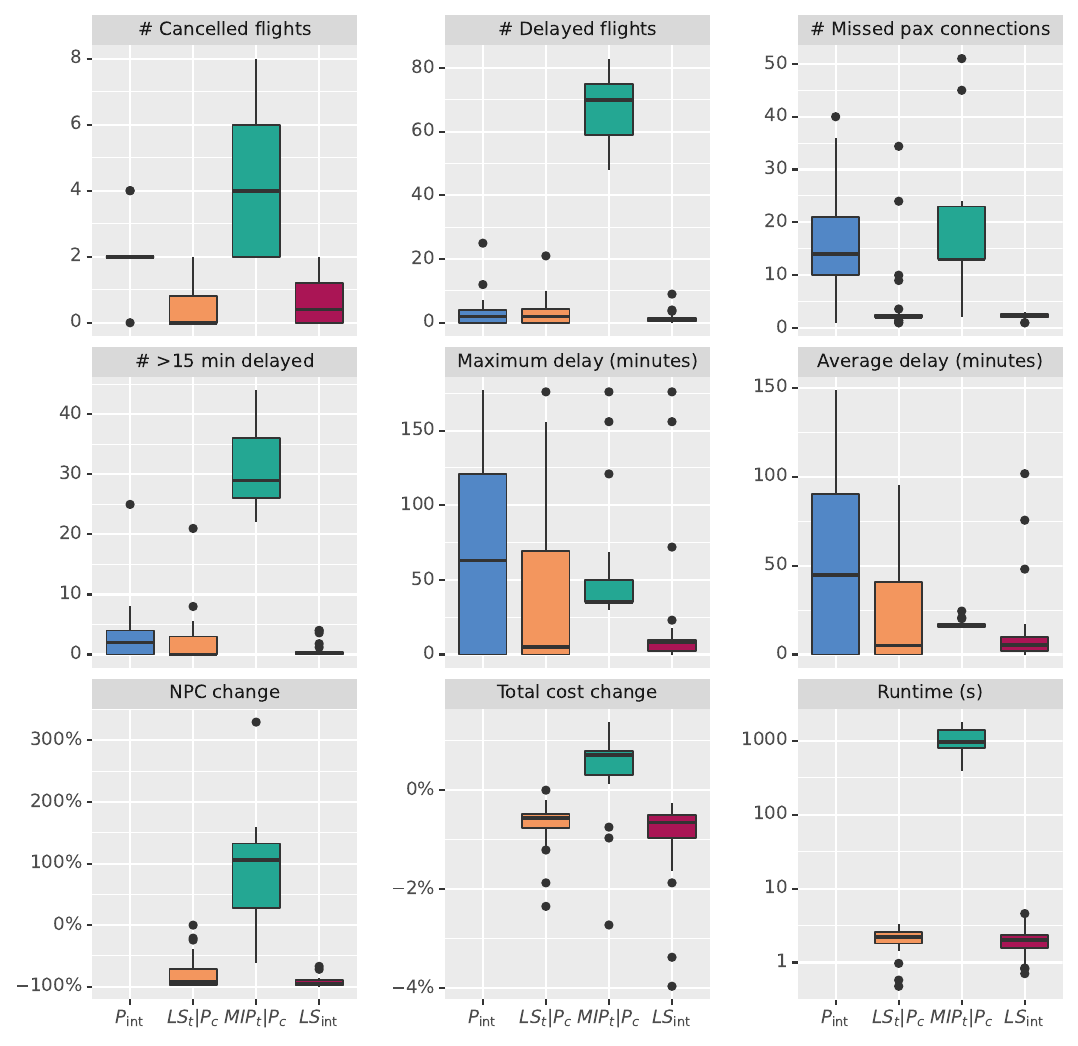}
    \caption{Comparison of KPI values when recovering both tail and crew schedules in case of the UR scenario.}
    \label{fig:kpi_boxplot_ur}
\end{figure*}

Next, we consider the integrated problem.
We look in more detail into the quality of the solutions.
This comparison is shown in \Cref{tab:solution_quality_kpi}.
Note that here, we show the average KPI values for disruptions that are solved by all four algorithms.
We also show how these values are distributed in \Cref{fig:kpi_boxplot_fd,fig:kpi_boxplot_ur}, showing these values in case of flight delays and unavailable resources respectively.
We first note that the local search methods can solve more instances than other approaches, 90-100$\%$ vs. 40-60$\%$.
The local search methods are also the only ones to improve cost compared to delay propagation.
Since the total cost encompasses the operating costs of three full days, these small percentage cost reductions result in big cost reductions.
It is crucial to note that the runtime of our integrated local search is below 30 seconds, which is feasible for real-time application.
Moreover, it leads to a great reduction of NPC, about 40$\%$ compared to delay propagation.
Even though on average $LS_\text{int}$ delays more flights, it has the least number of missed passenger connections and the lowest average delay.
Cost-wise, it is a lot better to take these connections into account, and, although counter-intuitive, solutions with a lower total delay may be more costly.

Another observation is that the $MIP_t|P_c$ model unexpectedly performs a lot worse than delay propagation.
This model returns the best possible tail recovery, however this clearly has a negative effect on the resulting integrated solution, adding significant delays to the solution.
Compared to $LS_t|P_c$ the $MIP_t|P_c$ model performs a lot worse, even though both models work on the same principle.
Looking into this, we see that $MIP_t|P_c$ performs four times as many tail assignment changes compared to $LS_t|P_c$.
This is likely because $LS_t$ does not look at parts of the schedule where no disruptions occur, it only looks at the issues present, where $MIP_t$ looks at the complete network and makes any change necessary in order to improve the tail recovery.
As shown in \Cref{tab:tail_quality_kpi}, this indeed decreases the total cost of a solution.
However, it comes at the cost of creating more issues for crew recovery.
To confirm this, we ran all methods with a penalty for making changes to the original assignment.
In \Cref{tab:solution_quality_kpi_with_swap_penalty} we compare the cost difference between all solution methods.
Setting a penalty on making assignment changes does indeed impact the $MIP_t$ model a lot as it now performs better than $LS_t|P_c$.
However, our $LS_\text{int}$ method still performs much better.
All in all, this shows that not incorporating any crew data clearly results in worse schedules, showing the benefit of an integrated solver.

\begin{table*}[ht!]
\centering
\caption{Comparison of average KPI values when recovering both tail and crew schedules.}
\label{tab:solution_quality_kpi}
\begin{tabular}{llrrrr}
\toprule
 &  & $P_{\text{int}}$ & $LS_{t}|P_{c}$ & $MIP_t|P_c$ & $LS_{\text{int}}$ \\
Disruption & KPI &  &  &  &  \\
\midrule
\multirow[c]{10}{*}{FD} & \# Cancelled flights & 0.00 & 0.10 & 4.64 & 0.18 \\
 & \# Delayed flights & 42.22 & 61.68 & 108.38 & 44.44 \\
 & \# Missed pax connections & 75.94 & 54.20 & 66.06 & 44.72 \\
 & \# $> 15$ min delayed & 19.78 & 25.23 & 52.84 & 13.18 \\
 & Maximum delay (minutes) & 54.72 & 57.31 & 67.16 & 46.84 \\
 & Average delay (minutes) & 18.05 & 16.46 & 18.55 & 13.62 \\
 & NPC change & n.a. & -18.25\% & +46.81\% & -41.61\% \\
 & Total cost change & n.a. & -0.44\% & +0.70\% & -1.22\% \\
 & Runtime (s) & n.a. & 14.85 & 1162.70 & 16.47 \\
 & Percentage solved & 64\% & 94\% & 64\% & 100\% \\
\cmidrule[0.01em]{1-6}
\multirow[c]{10}{*}{UR} & \# Cancelled flights & 2.07 & 0.57 & 4.07 & 0.55 \\
 & \# Delayed flights & 3.14 & 3.01 & 67.72 & 1.37 \\
 & \# Missed pax connections & 15.07 & 4.44 & 16.79 & 2.15 \\
 & \# $> 15$ min delayed & 2.86 & 2.06 & 30.93 & 0.61 \\
 & Maximum delay (minutes) & 68.21 & 36.20 & 50.07 & 19.86 \\
 & Average delay (minutes) & 49.00 & 21.79 & 16.77 & 12.51 \\
 & NPC change & n.a. & -79.98\% & +88.51\% & -92.83\% \\
 & Total cost change & n.a. & -0.70\% & +0.44\% & -0.95\% \\
 & Runtime (s) & n.a. & 2.12 & 1082.81 & 2.13 \\
 & Percentage solved & 43\% & 91\% & 43\% & 100\% \\
\bottomrule
\end{tabular}
\end{table*}

\begin{table}[ht!]
\centering
\caption{Comparison of average cost change when a penalty is applied for changing resource assignments.}
\label{tab:solution_quality_kpi_with_swap_penalty}
\begin{tabular}{llrrrr}
\toprule
 &  & $P_{\text{int}}$ & $LS_{t}|P_{c}$ & $MIP_t|P_c$ & $LS_{\text{int}}$ \\
Disruption & KPI &  &  &  &  \\
\midrule
\multirow[c]{2}{*}{FD} & NPC change & n.a. & -21.03\% & -3.03\% & -34.86\% \\
 & Total cost change & n.a. & -0.55\% & -0.64\% & -1.04\% \\
\cmidrule[0.01em]{1-6}
\multirow[c]{2}{*}{UR} & NPC change & n.a. & -78.34\% & -50.19\% & -89.59\% \\
 & Total cost change & n.a. & -0.66\% & -0.83\% & -0.91\% \\
\bottomrule
\end{tabular}
\end{table}

\begin{figure*}[ht!]
    \centering
    \begin{subfigure}[t]{0.48\textwidth}
        \centering
        \includegraphics[width=\linewidth]{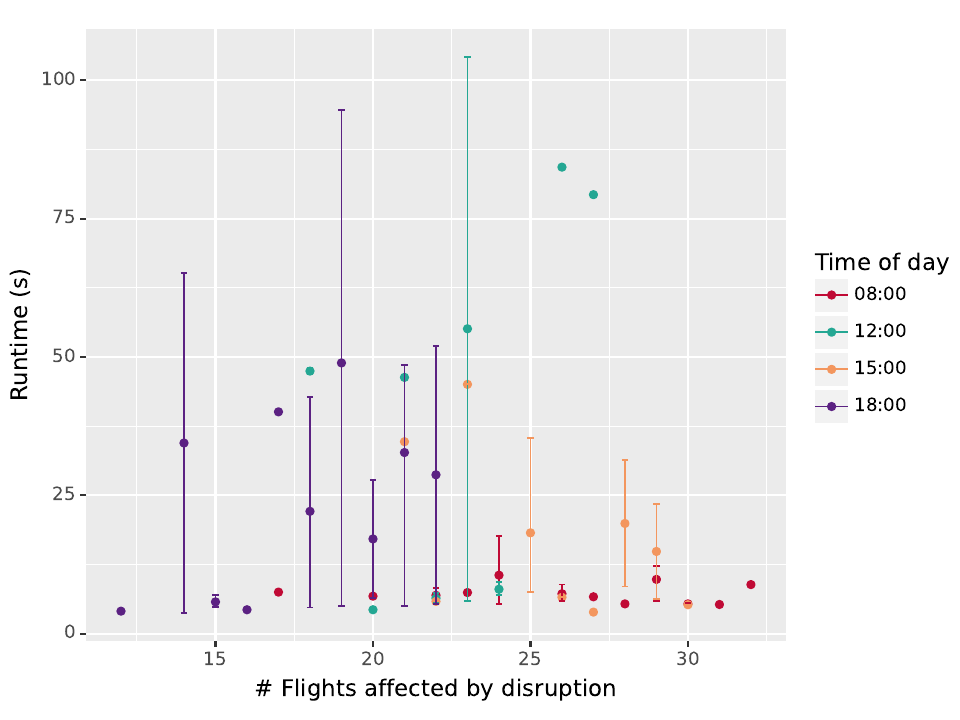}
        \caption{Average runtime of the simulated annealing compared to the size of the disruption. Next to that the $95\%$ confidence interval of this average is displayed.}
        \label{fig:disruption_correlation_runtime}
    \end{subfigure}
    \hfill
    \begin{subfigure}[t]{0.48\textwidth}
        \centering
        \includegraphics[width=\linewidth]{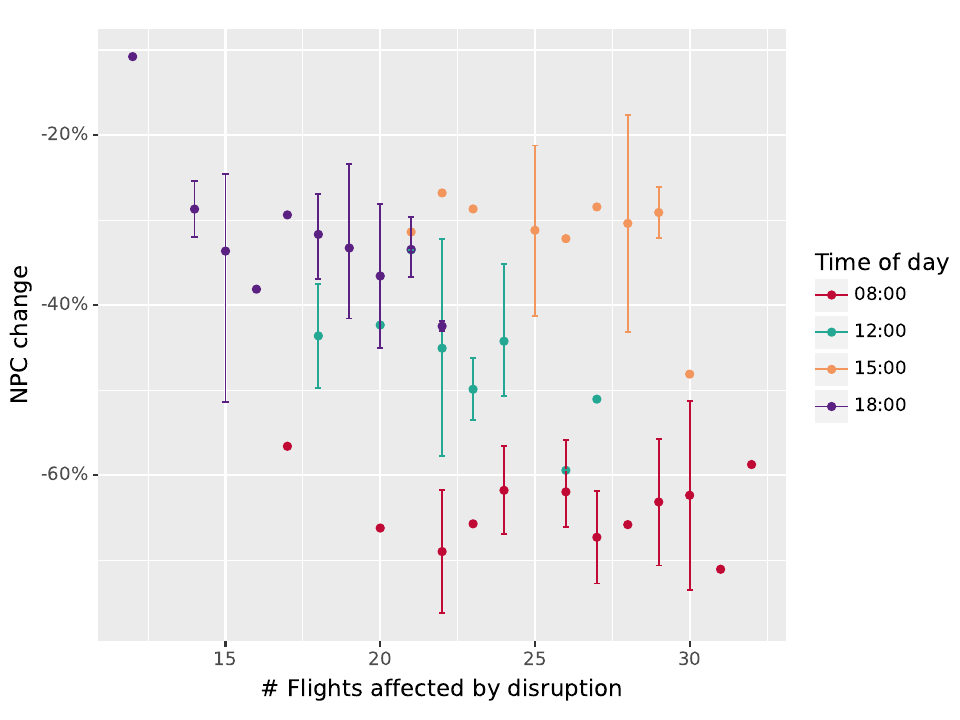}
        \caption{Average NPC reduction compared to $P_\text{int}$ compared to the size of the disruption. Next to that the $95\%$ confidence interval of this average is displayed.}
        \label{fig:disruption_correlation_npc}
    \end{subfigure}
    \caption{Comparison between the size of a disruption, measured in the number of flights that are affected, and either the runtime of the simulated annealing or the reduction in NPC compared to $P_\text{int}$.}
    \label{fig:disruption_correlation}
\end{figure*}

Apart from good solution quality, it is also important to know how the algorithm reacts under different disruption sizes.
While we do not generate relatively big and relatively small disruptions separately, our disruption generation has some variation in the size of disruption it generates.
For the size of the disruption, we measure the number of flights that are affected by the disruption.
These flights need a decision as they otherwise could not be flown.
We compare this size to the runtime of our algorithm and to the NPC change compared to $P_\text{int}$.
The comparisons are shown in \Cref{fig:disruption_correlation_runtime,fig:disruption_correlation_npc} respectively.
In here, we first notice the runtime to behave quite erratic, and we do not see a correlation between the runtime and the size of the disruption.
The NPC change stays somewhat similar given the time of the disruption.
Meaning that given the time of the disruption a similar (percentage) improvement can be made by our algorithm compared to $P_\text{int}$.
Note that a bigger disruption will result in higher NPC costs, but the percentage improvement stays similar.
From these results we see that our algorithm shows similar performance across different disruption sizes, although more elaborate testing may be done to confirm this.

\subsection{Comparing Different Cost Models}\label{sec:results_cost_model}
Next, we investigate the effects of different choices within the recovery model.
Here, we mainly investigate how our algorithm performs under different cost models and how it can be used to fine-tune the recovery approach.
Different airlines might have different priorities on how they want to recover their schedules.
As can be seen in the previous results, the current cost model sometimes favours saving passenger connections over delays, and hence over schedule stability.
In this section, we look at two priorities an airline might have for how to recover their schedules.
These are given by the cost model for delays, and a focus on reducing the number of reassignments in the final schedule.

To get a realistic approximation of the various reaccommodation and compensations costs accumulated when delaying and cancelling flights, our model uses precomputed delay costs and a cost for breaking passenger connections.
Not only does this provide a realistic approximation, these costs also express passenger dissatisfaction in the model, and are also a way to prioritize certain flights and connections over others.
However, in a recovery scenario, the inclusion of these costs might result in more flights being delayed than necessary, meaning that we might compromise on operational stability.
Here, the airline has a choice in how they want to prioritize their recovery.
For this reason, we compare different ways for including delay costs in our model.
We compare the original approach, that uses these precomputed delay costs, with a more punctuality-driven approach.
In this punctuality driven approach, we penalize any flight if it has a delay.
However, flights with a delay greater than \num{15} minutes are penalized more, as these flights are usually not considered punctual any more.
For example, \citet{PRREuroControl2024} uses a \num{15}-minute threshold to determine punctuality.
Our final delay cost function uses a constant penalty for a delayed flight and this penalty is multiplied by 10 if the delay is greater than \num{15} minutes.
Furthermore, passenger connections are not considered.
To compare cost, we evaluate the results of both cost models using the original cost model, which was also used in our previous results.
This cost model should represent a more realistic picture of the actual costs incurred.

\begin{figure*}[ht!]
    \centering
    \includegraphics[width=0.85\textwidth]{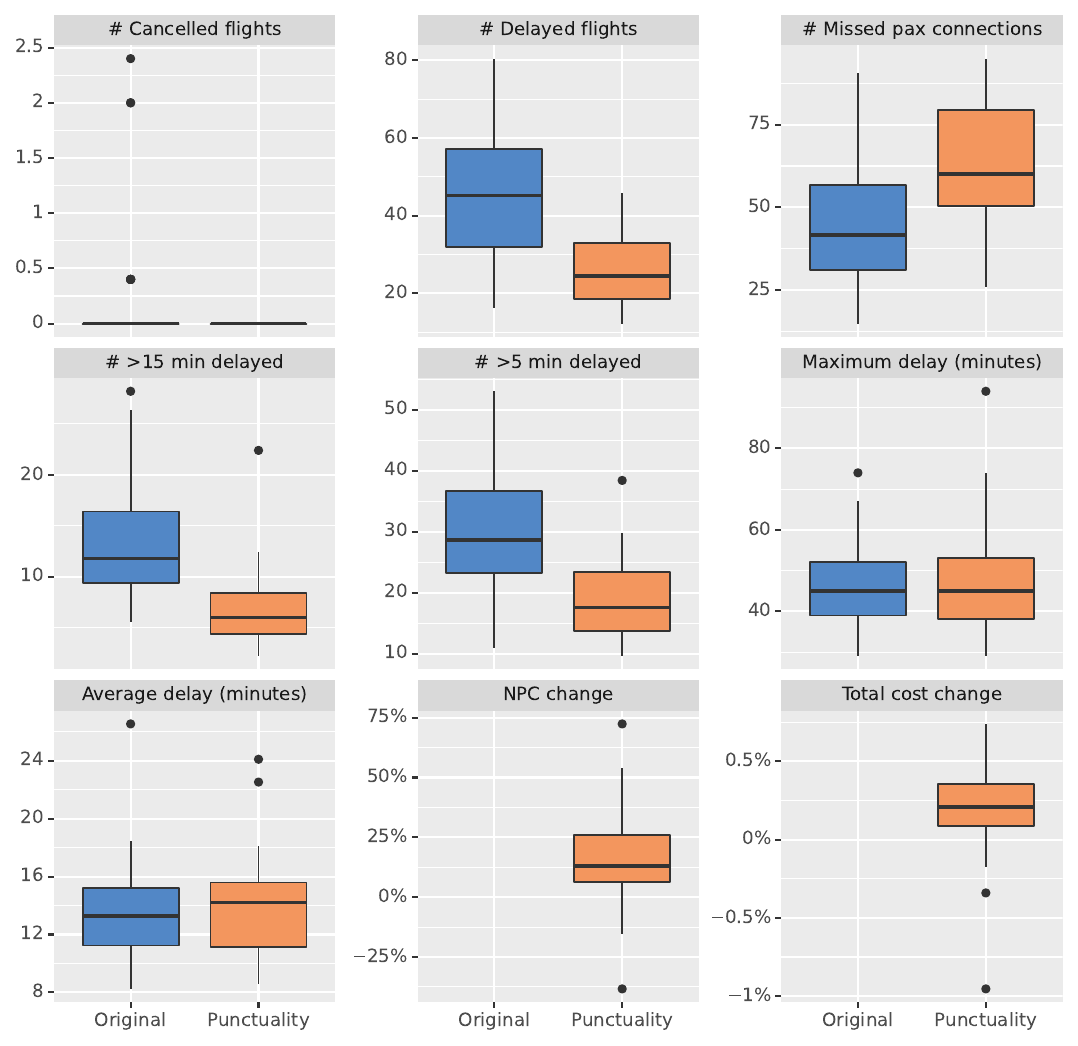}
    \caption{Comparison of KPI values for various delay cost models in case of the FD scenario.}
    \label{fig:cost_model_kpi_boxplot_fd}
\end{figure*}

\begin{figure*}[ht!]
    \centering
    \includegraphics[width=0.85\textwidth]{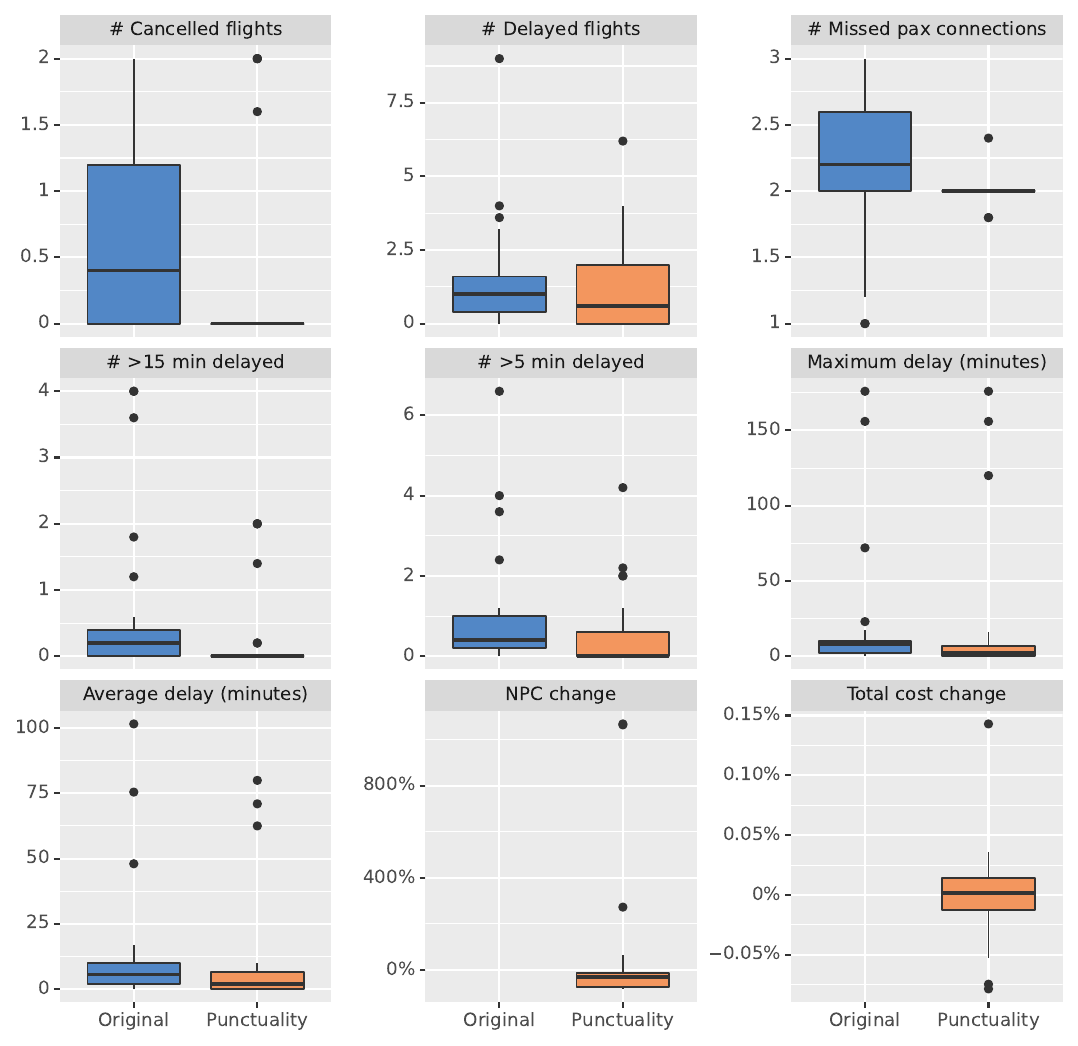}
    \caption{Comparison of KPI values for various delay cost models in case of the UR scenario.}
    \label{fig:cost_model_kpi_boxplot_ur}
\end{figure*}

The comparison between average KPI values is shown in \Cref{tab:punctuality_cost_comparison}.
Furthermore, we show these values are distributed in \Cref{fig:cost_model_kpi_boxplot_fd,fig:cost_model_kpi_boxplot_ur}, for the FD and UR disruptions respectively.
These results show that the use of punctuality costs indeed leads to fewer and shorter delays, and also reduces the number of cancellations.
However, neglecting the passenger connections leads, in this case, to a higher NPC and overall costs.
In case of unavailable resources, the results are not as conclusive, with both models returning the same total cost.
The results seem near-optimal, as the delays are quite small.
We see a difference in the KPIs, where the original cost model prefers a cancellation over less delay, while for punctuality costs it is the other way around.

\begin{table*}[ht!]
\centering
\caption{Comparison of different KPIs for various delay cost models.}
\label{tab:punctuality_cost_comparison}
\begin{tabular}{llrr}
\toprule
 &  & Original & Punctuality \\
Disruption & KPI &  &  \\
\midrule
\multirow[c]{9}{*}{FD} & \# Cancelled flights & 0.19 & 0.00 \\
 & \# Delayed flights & 44.84 & 26.62 \\
 & \# $> 15$ min delayed & 13.27 & 6.90 \\
 & \# $> 5$ min delayed & 29.83 & 18.95 \\
 & Maximum delay (minutes) & 46.24 & 46.64 \\
 & Average delay (minutes) & 13.51 & 13.93 \\
 & \# Missed pax connections & 44.66 & 62.30 \\
 & NPC change & n.a. & +15.88\% \\
 & Total cost change & n.a. & +0.21\% \\
\cmidrule[0.01em]{1-4}
\multirow[c]{9}{*}{UR} & \# Cancelled flights & 0.55 & 0.33 \\
 & \# Delayed flights & 1.37 & 1.18 \\
 & \# $> 15$ min delayed & 0.61 & 0.27 \\
 & \# $> 5$ min delayed & 0.99 & 0.53 \\
 & Maximum delay (minutes) & 19.86 & 18.53 \\
 & Average delay (minutes) & 12.51 & 9.66 \\
 & \# Missed pax connections & 2.15 & 2.00 \\
 & NPC change & n.a. & +50.34\% \\
 & Total cost change & n.a. & 0.00\% \\
\bottomrule
\end{tabular}
\end{table*}

\begin{figure*}[ht!]
    \centering
    \includegraphics[width=0.7\linewidth]{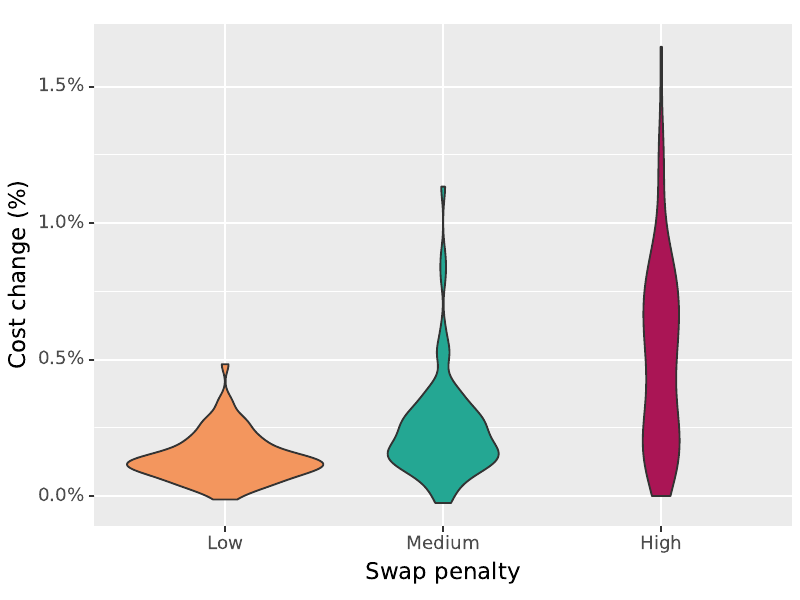}
    \caption{Comparison between different swap penalty severities, showing the cost change compared to not using a swap penalty.}
    \label{fig:swap_cost_cost_change}
\end{figure*}

\begin{figure*}[ht!]
    \centering
    \includegraphics[width=0.7\linewidth]{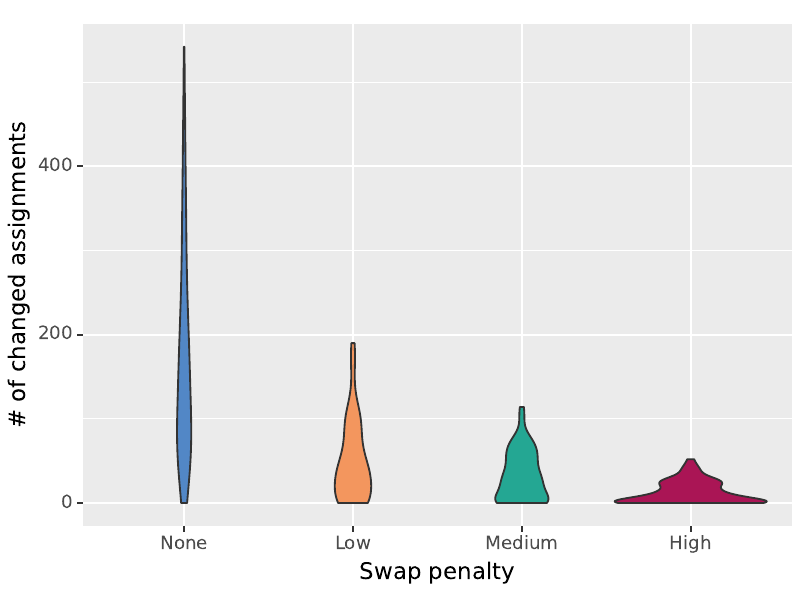}
    \caption{Comparison between different swap penalty severities, showing the number of changed assignments in the final schedule.}
    \label{fig:swap_cost_number_changed_assingments}
\end{figure*}

An important consideration in disruption management is the number of assignment changes proposed to recover a schedule.
It is often preferred to keep the schedule as close to the original schedule as possible.
We limit the number of swaps by imposing a penalty on the number of changed assignments, which we call the ``swap penalty'', and investigate the effect this has on the solution cost.
To investigate this, we only look at the FD scenarios, as these seem to be more interesting.
First, we compare cost changes when changing the swap penalty in \Cref{fig:swap_cost_cost_change}.
Here we show the percentage change compared to the case without a swap penalty.
From this figure, we see that introducing such a penalty will indeed lead to higher costs.
Next, we compare the number of changed assignments made in the final solution.
This is shown in \Cref{fig:swap_cost_number_changed_assingments}.
From this figure, it is clear that varying the severity of the swap penalty has a significant impact on the total number of changed assignments.
This does increase the cost, as we saw previously, but this might be worth the price.

\subsection{Impact of type swaps}\label{sec:results_type_swap}
As mentioned in \Cref{sec:local_search_neighbours_swaps}, our algorithm allows for so-called type swaps, i.e. swapping an aircraft with one of a different type.
This creates extra flexibility in the recovery, especially for airlines that fly different aircraft types of similar sizes.
When changing to an aircraft of different size many passenger-related considerations must be taken into account, as either some may be stranded or the aircraft itself has more empty seats.
In our data, KLM does not yet have aircraft of different types with a similar size.
However, it will be more relevant for them in the future when they phase out the Boeing 737 for the Airbus A320.
Still, in big disruptions this may be an important lever to have.
Thus, we compare how our algorithm reacts in case we allow or disallow these type swaps.

A comparison of different KPIs in both these scenarios is given in \Cref{tab:multi_fleet_comparison}.
In our final solutions, we see that on average there are about one or two of these swaps present.
Furthermore, we see that in the FD scenario the addition of type swaps indeed leads to slightly better and more cost-efficient solutions.
The UR scenario is very similar to what we have seen earlier in \Cref{sec:results_neigbhour_generation_strategies,sec:results_cost_model}, where a different solution can lead to very big changes in the NPC with the overall cost being very similar.
Overall we see promising benefits for considering these type swaps, especially in the case of many delayed flights.
In our current data, we do not have mixed fleet of aircraft, although these mixed-fleet operations are coming.
With this consideration, our results show potential for airlines operating a mixed-fleet, where these type swaps are likely to be more beneficial.

\begin{table}[ht!]
\centering
\caption{Comparison of different KPIs for allowing or disallowing type swaps.}
\label{tab:multi_fleet_comparison}
\begin{tabular}{llrr}
\toprule
 &  & No Type Swaps & Allow Type Swaps \\
Disruption & KPI &  &  \\
\midrule
\multirow[c]{9}{*}{FD} & \# Cancelled flights & 0.15 & 0.14 \\
 & \# Delayed flights & 50.61 & 51.32 \\
 & \# Missed pax connections & 52.41 & 51.71 \\
 & \# $> 15$ min delayed & 14.78 & 14.57 \\
 & Maximum delay (minutes) & 52.13 & 52.11 \\
 & Average delay (minutes) & 13.60 & 13.45 \\
 & NPC change & n.a. & -1.39\% \\
 & Total cost change & n.a. & -0.03\% \\
 & Runtime (s) & 19.77 & 19.77 \\
\cmidrule[0.01em]{1-4}
\multirow[c]{9}{*}{UR} & \# Cancelled flights & 0.37 & 0.41 \\
 & \# Delayed flights & 1.62 & 1.93 \\
 & \# Missed pax connections & 2.32 & 2.38 \\
 & \# $> 15$ min delayed & 0.82 & 0.83 \\
 & Maximum delay (minutes) & 26.91 & 27.57 \\
 & Average delay (minutes) & 19.07 & 17.83 \\
 & NPC change & n.a. & +16.83\% \\
 & Total cost change & n.a. & +0.01\% \\
 & Runtime (s) & 2.20 & 2.60 \\
\bottomrule
\end{tabular}
\end{table}

%% file: 6_conclusion.tex
In this paper, we present a simulated annealing algorithm to solve the integrated airline fleet and crew recovery problem with considerations for passenger impact.
In our experiments, the algorithm has an average runtime below 30 seconds, and a median runtime below 10 seconds.
Together with the high degree of parallelizability, this makes our algorithm suitable for time-critical operational use.
Our simulated annealing algorithm uses directed search, in the sense that each iteration of the algorithm involves an issue in the current schedule.
We presented different strategies to generate neighbours, where our candidate selection methods performed generally better in the flight delay scenarios, and a random neighbour selection is a bit better in case of unavailable resources.

To evaluate the quality of our solutions, we first compare local search for tail recovery with a MIP model.
We show that local search finds near-optimal solutions while running ten times faster.
We then compare our complete algorithm with two sequential approaches: one that uses local search for tail recovery and delay propagation for crew recovery, and one that uses a MIP model for tail recovery and delay propagation for crew recovery.
Our experiments suggest that the integrated approach reduces total delay and costs, while also solving more cases.
Compared to solving disruptions in an integrated way, but by using only delay propagation, it reduces missed passenger connections and average delay, resulting in a 40\% reduction in non-performance costs.
These results show the benefit of considering crew resources during tail recovery.

We also compare cost models that reflect different airline priorities during schedule recovery: we compare our original model, which accounts for passenger connections with a model that prioritizes minimizing delays and with one that prioritizes minimizing the number of reassignments.
Our experiments show that across these cost models, the resulting solutions reflect the different priorities encoded in each cost model, as expected.
These comparisons can help airlines to determine a cost model that aligns with their priorities.

The practical relevance of this work is also reflected in feedback from Bob Tulleken, Vice President of Operations and Decision Support at KLM:
\begin{quote}
The integrated recovery algorithm presented in this paper is a possible game-changer for KLM.
By addressing disruptions in both fleet and crew schedules simultaneously, it offers a holistic approach to minimizing operational costs and passenger impact.
The fast runtime ensures that decisions can be made in real-time, enhancing our ability to maintain schedule integrity and customer satisfaction.
This innovative solution not only reduces Non-Performance Costs but also exemplifies the benefits of integrating multiple resources into our planning processes, promising highly valuable results for KLM.
\end{quote}
The quote highlights how integrating fleet and crew recovery, together with a short runtime, can support real-time decision-making.

In our experiments, we test on a flight network with a single hub and short- and medium-haul flights.
Although we do not explicitly evaluate multiple-hub networks, our neighbourhood structures are designed to support them.
Extending the model to multiple hubs would require accounting for resource balancing between the hubs.
To extend the model to long-haul flights, a larger recovery window is required, since long-haul flights take longer.
This may lead to an increase in runtime; however, the magnitude of this effect should be assessed in future experiments.
In case of long-haul flights, crew connections become less important, while the location of crew becomes more important.
The model would therefore need to be extended with crew-repositioning options, such as deadheading.
Together with the flexibility in the cost model, these extensions suggest that our model could be adapted to a broad range of airline operations.

Given the runtime and solution quality of our local search algorithm, several extensions merit further research.
First, given the importance of crew positioning in long-haul operations, future work could model crew deadheading explicitly.
This could not only support long-haul operations but also provide additional recovery options during major disruptions.
Further integration could include maintenance assignments, which are currently treated as fixed in our model.
In practice, these assignments may sometimes be changed, depending on safety considerations and the availability of technicians and spare parts.
Another extension would be to include cabin crew, while accounting for their specific rules which are slightly different from cockpit crew.
Further work could also integrate ground resources or add robustness measures, allowing the algorithm to not only resolve disruptions but also reduce the impact of future ones.

%% file: 7_tail_mip.tex
We introduce a MIP formulation that can be used for solving disruptions in only the tail schedule.
We use this formulation when benchmarking our integrated algorithm against sequential approaches.
For this, several different formulations based on different network models can be used \citep{Clausen2010,Hassan2021}.
In our formulation we use a connection network similar to the work of \citet{Andersson2004}, \citet{Liang2018}, and others.

In a connection network, we consider a directed graph $G = (V, A)$, where $V$ consists of all flight activities $F$, maintenance activities $M$, and reserve activities.
Thus, each $v \in V$ represents an activity to be performed.
Then an arc $(v_1, v_2) \in A$ exists if and only if activities $v_1$ and $v_2$ can be performed successively by the same resource.
This problem can now be modelled as a flow problem, where the schedule for each aircraft is represented by a unit flow.
Therefore, the graph $G$ is extended with source and sink nodes $s$ and $t$.

For this formulation, we introduce the following variables:
\begin{align*}
    x_{a,r} &\in \{0, 1\} && \text{Denotes if tail $r \in T$ is assigned to arc $a \in A$,} \\
    y_{f} &\in \{0, 1\} && \text{Denotes if flight $f$ is cancelled,} \\
    d_{f} &\in [0, 180] && \text{The delay of flight $f$,} \\
    z_{p} &\in \{0, 1\} && \text{Denotes whether passenger connection $p$ is broken.}
\end{align*}

The objective is the same as \Cref{eq:objective}.
However, it needs to be adapted to the introduced assignment variables.
Thus, the objective is to minimize
\begin{equation}
    \sum_{(f',f) \in A}\sum_{r \in T} x_{f'f,r} \left( c_{f,r}^{\text{op}} + c^{\text{swap}}_{f,r} \right)
    + \sum_{f \in F} y_{f} c_{f}^{\text{cancel}}
    + \sum_{f \in F} c_{f}^{\text{delay}}(d_{f})
    + \sum_{p \in P} z_{p} c_{p}^{\text{mc}}
.\end{equation}

This is subject to the constraints listed below.
Since the problem is modelled as a flow problem, first the flow constraints are added.
These are the constraints:
\begin{align}
    \sum_{\substack{f' \in V,\\(f', f) \in A}} x_{f'f,r} &= \sum_{\substack{f' \in V,\\(f, f') \in A}} x_{ff',r} & \forall f \in V, r \in T, \\
    \sum_{\substack{f \in V,\\(s, f) \in A}} x_{sf,r} &= 1 & \forall r \in T, \\
    \sum_{\substack{f \in V,\\(f, t) \in A}} x_{ft,r} &= 1 & \forall r \in T.
\end{align}
These constraints handle the flow preservation and ensure that each aircraft corresponds to a single flow.

\paragraph{Flight Assignment}
A flight needs to be either assigned or cancelled, which is done with the constraint
\begin{equation}
    y_{f} + \sum_{\substack{f' \in V,\\ (f', f) \in A}}\sum_{r \in T} x_{f'f,r} = 1
    \qquad\qquad
    \forall f \in V
.\end{equation}
Here we note that for a maintenance activity $m \in M$ that needs to be performed on tail $t \in T$, we have that $y_m = 0$, enforcing that this activity cannot be cancelled.
Furthermore, we need to add an additional constraint to enforce that the activity is performed with tail $r$.
This is the constraint
\begin{equation}
     \sum_{\substack{f' \in V,\\ (f', m) \in A}} x_{f'm,r} = 1
     \qquad\qquad\qquad
     \forall m \in M
.\end{equation}

\paragraph{Flight Delays}
To set the correct flight delay and to ensure that the turn-around time for aircraft is respected, we make use of big-M constraints.
The turn-around time between two flights is only dependent on the subtypes, which are specific aircraft variants.
Thus, let $S$ be the set of subtypes, and $t^{\text{conn}}_{ff',s}$ the turn-around time between flights $f$ and $f'$ for subtype $s$.
Next, we define $T_s = \{ r \in T \mid r \text{ of subtype } s \}$, meaning that $T_s$ contains all tails of the specified subtype $s$.
Then to enforce correct delays and maintain the turn-around time, we add the big-M constraint
\begin{equation}
    (f'_d + d_{f'}) - (f_a + d_{f}) \geq t^{\text{conn}}_{ff',s} - M \left( 1 - \sum_{r \in T_s} x_{ff',r} \right)
    \qquad
    \forall (f, f') \in A, s \in S
.\end{equation}

\paragraph{Missed Passenger Connections}
To detect that passenger connections are broken, we need another big-M constraint.
Here we use $t^{\text{pax}}_{ff'}$ to denote the minimum connection time between the flights $f$ and $f'$, where $(f, f') \in P$.
This minimum connection time can be broken due to delays, or when one of the flights is cancelled.
Thus, we include constraints to set a broken passenger connection when one of the flights is cancelled or when the minimum connection time does not hold.
This means that the following constraints are added to the model:
\begin{align}
    (f'_d + d_{f'}) - (f_a + d_{f}) &\geq t^{\text{pax}}_{ff'} - M z_{ff'} && \forall (f, f') \in P, \\
    z_{ff'} &\geq y_{f}  && \forall (f, f') \in P, \\
    z_{ff'} &\geq y_{f'} && \forall (f, f') \in P.
\end{align}